%% file: main.tex
\documentclass[hidelinks,onefignum,onetabnum]{siamart250211}

\usepackage[utf8]{inputenc}
\usepackage[T1]{fontenc}
\usepackage{lmodern}
\usepackage[english]{babel}

\usepackage{tumcolors}
\usepackage{tummacros}
\colorlet{blue}{tumblue}

\usepackage{booktabs} 
\usepackage{graphicx} 
\usepackage{pgfplots} 
\usepackage{amssymb, mathtools} %
\pgfplotsset{compat=1.18}
\usepackage{hyperref} 
\usepackage{cleveref} 
\usepackage{comment}

\usepackage{tabularx}

\usepackage{subcaption}
\usepackage{caption}

\usepackage[obeyFinal, backgroundcolor=white, linecolor=tumgray, textcolor=red, textsize=footnotesize]{todonotes}

\providecommand{\data}[1]{{\par\small\noindent\textbf{\textit{Data availability:}} #1}}
\providecommand{\funding}[1]{{\par\small\noindent\textbf{\textit{Funding:}} #1}}

\newcommand\titlename{A kernel compression method for distributed-order fractional partial differential equations}

\hypersetup{
    pdftitle={\titlename},
	pdfauthor={Jonas Beddrich, Barbara  Wohlmuth},
	colorlinks=true,
	linkcolor=\linkcolor,
	urlcolor=\urlcolor,
	citecolor=\citecolor,
	menucolor=\menucolor
}

\hypersetup{
  pdftitle={A kernel compression method for distributed-order fractional partial differential equations},
  pdfauthor={J. Beddrich and B. Wohlmuth}
}

\author{Jonas Beddrich \thanks{Technical University of Munich, School of Computation, Information and Technology, Department of Mathematics, Boltzmannstra{\ss}e 3, D-85748, Garching b. München, Germany (\email{beddrich@ma.tum.de}).}
\and Barbara Wohlmuth \footnotemark[2]
}
\title{A kernel compression method for distributed-order fractional partial differential equations \thanks{Submitted on 18 August 2025.
\funding{Nothing to declare.}}}

\begin{document}

\maketitle

\begin{abstract}
    \input{Files/abstract}
\end{abstract}



\data{A GitHub link for reproducibility will be provided after acceptance of the manuscript.}


\section{Introduction}
\label{sec:intro}
\input{Files/Introduction}

\section{Kernel compression for distributed-order differential operators}
\label{sec:kernelcompression}
\input{Files/KernelCompression}

\section{Exponential sum approximation}
\label{sec:exponentialsumapproximation}
\input{Files/ExponentialSums}

\section{Numerical method}
\label{sec:numerical_scheme}
\input{Files/NumericalMethod}

\section{Numerical results}
\label{sec:numerical_experiments}
\input{Files/NumericalResults}

\section{Conclusion}
\label{sec:conclusion}
\input{Files/Conclusion}

\bibliographystyle{siamplain}
\bibliography{references}

\end{document}

%% file: Files/abstract.tex
We propose a kernel compression method for solving Distributed-Order (DO) Fractional Partial Differential Equations (DOFPDEs) at the cost of solving corresponding local-in-time PDEs. 
The key concepts are (1) discretization of the integral over the order of the fractional derivative and (2) approximation of linear combinations of integral kernels with exponential sums, expressing the non-local history term as a sum of auxiliary variables that solve a weakly coupled, local in time system of PDEs. 
For the second step, we introduce an improved algorithm that approximates the occurring integral kernels with double precision accuracy using only a moderate number ($<100$) of exponential terms. 
After temporal discretization using implicit Runge--Kutta methods, we exploit the inherent structure of the PDE system to obtain the solution at each time step by solving a single PDE. 
At the same time, the auxiliary variables are computed by a linear update, not even requiring a matrix-vector multiplication. 
Choosing temporal meshes with a grading factor corresponding to the convergence order of the Runge--Kutta schemes, we achieve the optimal decay of the temporal discretization error.
The flexibility and robustness of our numerical scheme are illustrated by recreating well-studied test cases and solving linear and nonlinear DOFPDEs in 2D and 3D with up to 40 million spatial degrees of freedom. 

%% file: Files/Introduction.tex
Time-fractional differential equations (FDEs) play a crucial role in modeling anomalous transport and damping effects. 
However, realistic dynamics often involve heterogeneous materials and the interaction of multiple physical processes, which are poorly described using a single time-fractional derivative.  
To capture these complex mechanisms, the more general concept of DO Fractional Diffferential Equations (DOFDEs) was introduced, whereby the DO differential operator is defined as a weighted mean of time-fractional derivatives \cite{caputo1966linear,caputo1967linear}.  
In recent years, DOFDEs have been successfully applied to model groundwater flow \cite{yin2020distributed}, oxygen delivery through capillary tissues \cite{srivastava2010multi}, dental restoration materials \cite{petrovic2015viscoelastic}, and to characterize relaxation processes in magnetic resonance imaging \cite{yu2022application}, motivating the development of efficient numerical methods for their simulation. 

We consider DO differential operators with weight function $\phi(\alpha)$ of the form 
\begin{equation}
    \label{eq:IN_DODE}
    \int_0^\alphamax \phi(\alpha) D^\alpha_C u(t) \Diff \alpha, 
\end{equation}
where $\alphamax \in \N$ is the smallest natural number such that $\mathrm{supp}(\phi) \subset [0,\alphamax]$, and $D^\alpha_C$ is the time-fractional derivative of Caputo-type. 
For $n \in \N$, $D^n_C f(t) := f^{(n)}(t) := \frac{\Diff^n}{\Diff t^n}f(t) $, while for a non-integer $\alpha > 0$, it is defined as  
\begin{equation}
    \label{eq:def_Caputo}
    D^\alpha_C f(t) 
    := g_{\lceil\alpha\rceil-\alpha} * f^{(\lceil\alpha\rceil)} (t). 
\end{equation} 
Here, $g_\alpha(t) = t^{\alpha-1} / \Gamma(\alpha)$ is the Riemann--Liouville (RL) fractional integral kernel and $\Gamma(\alpha) = \int_0^\infty x^{\alpha-1} e^{-x} \Diff x$ denotes the Gamma function \cite{diethelm2010analysis}. 
The DO differential operator \eqref{eq:IN_DODE} is a generalization covering constant-order time-fractional derivatives as well as weighted multi-term time-fractional differential operators, which we recover choosing $\phi(\alpha) = \delta(\alpha-\alpha_0)$ and $\phi = \sum_{n=1}^N \beta_n \delta(\alpha - \alpha_n)$, with fractional-order $\alpha_n > 0$ and weights $\beta_n$, respectively. Here, $\delta$ denotes the Dirac delta distribution. 

Consequently, the same challenges as for FDEs arise when solving DOFDEs. 
In particular, the non-locality of the DO differential operator necessitates efficient treatment of the memory term to overcome the associated significant computational cost and memory requirements of standard schemes. 
Further, independent of the regularity of the right-hand side, the solution of a DOFDE is possibly non-smooth at the initial time, requiring significant modification of spectral methods and numerical schemes depending on Taylor series expansions to achieve high-order accuracy. 
Instead of summarizing the vast literature addressing these problems for FDEs, we refer the interested reader to the two recent review articles \cite{diethelm2021new,diethelm2020good}. 

Numerical treatment of DO differential operators typically involves the numerical quadrature of the integral over the fractional order $\alpha$ and discretizing the resulting time-fractional derivatives, whose number matches the number of quadrature points. 
Various numerical methods have been proposed to solve the multi-term FDEs. 
By transformation into a high-dimensional system of FDEs, standard multi-step or Adams-type predictor-corrector methods can be applied \cite{diethelm2004multi,diethelm2009numerical}. 
A similar technique arises by treating all but one time-fractional derivative as fictitious source terms \cite{katsikadelis2014numerical}.
Numerous mesh-based schemes, such as the BDF method \cite{morgado2015numerical}, the L1 method \cite{bu2017finite}, or the Grünwald--Letnikov method \cite{aminikhah2019numerical}, as well as mesh-free techniques using Legendre \cite{yuttanan2019legendre}, Taylor wavelets \cite{yuttanan2021numerical}, or hybrid functions based on block-pulse functions in combination with Bernoulli \cite{mashayekhi2016numerical} and Taylor \cite{jibenja2018efficient} polynomials, have been applied to directly solve the multi-term FDE, and we refer to Section 2.4.2. of \cite{ding2021applications}, for an exhaustive overview of numerical methods applied in the context of DO differential operators. 
None of the methods mentioned above is well-suited for DOFPDEs, as both the non-locality and the large number of fractional derivative terms have to be managed efficiently. 

Our focus lies on kernel compression schemes, a group of non-classical methods that are both fast and memory efficient \cite{diethelm2022trends}.
The applicability of these methods to FPDEs has been demonstrated in the literature, e.g., \cite{beddrich2025numerical,beddrich2024numerical,khristenko2023solving}. 
The key concept is to replace the RL integral kernel with a sum of exponential terms
\begin{equation}
    \label{eq:sum_of_exponentials}
    g_{1-\alpha}(t) 
    = \frac{t^{-\alpha}}{\Gamma(1-\alpha)} 
    \approx \sum_{j=1}^m w_j e^{-\lambda_j t},
\end{equation}
where the weights $w_j$ and poles $\lambda_j$ are positive. 
By doing so, the non-local time-fractional derivative term is approximated by a sum of auxiliary variables given by a system of local in time PDEs.

To obtain such an exponential sum approximation, Yuan and Agrawal \cite{yuan2002numerical} substituted $x = ts^2$ in the definition of the gamma function, resulting in the integral identity 
\begin{equation}
    \label{eq:infinite_state}
     \frac{t^{-\alpha}}{\Gamma(1-\alpha)} = \frac{2}{\Gamma(\alpha)\Gamma(1-\alpha)} \int_0^\infty s^{2\alpha-1} e^{-ts^2} \Diff s, \quad 
     0 < \alpha < 1. 
\end{equation}
Subsequent discretization using the Gauss--Laguerre quadrature yields an approximation of $g_{1-\alpha}(t) $ in the form of \eqref{eq:sum_of_exponentials}, see also \cite{diethelm2008investigation}. 
This initial, suboptimal approach \cite{schmidt2006critique} was refined by transforming the integration domain to $(-1,1)$ and applying the Gauss--Jacobi quadrature \cite{diethelm2009improvement}. 
Many alternative changes of the integration variable \cite{beylkin2010approximation,chatterjee2005statistical,diethelm2021new,li2010fast,mclean2018exponential,singh2006galerkin} and quadrature schemes \cite{baffet2019gauss,banjai2019efficient,beylkin2005approximation,jiang2017fast,li2010fast,mclean2018exponential} have been proposed in the literature. 

More general methods that do not rely on integral identities, such as \eqref{eq:infinite_state}, and thus, the explicit form of the integral kernel, arise based on multipole approximations of the Laplace transformed integral kernel. 
Then, by taking the inverse Laplace transform, one recovers an exponential sum approximation in the form of \eqref{eq:sum_of_exponentials}. 
The multipole approximations have been computed based on complex contour integrals \cite{baffet2017high,lopez2008adaptive,lubich2002fast,schadle2006fast} or by rational approximation and subsequent partial fractions decomposition \cite{beddrich2025numerical,khristenko2023solving}, yielding good approximations with only a small number of exponential terms \cite{khristenko2023solving}. 
After application of a quadrature formula to a DO differential operator, naively approximating each resulting fractional derivative term using a kernel compression scheme is infeasible due to their large number and the subsequent amount of auxiliary variables. 
In this work, we overcome this challenge, introducing a kernel compression approach that approximates the DO differential operator using only a moderate number of auxiliary variables. 

The structure of the remaining article is as follows: 
In Section \ref{sec:kernelcompression}, we introduce the concept of kernel compression methods for time-fractional derivatives and generalize the concept to DO differential operators. 
An improved method for calculating the exponential sum approximation of the occurring integral kernels is presented in Section \ref{sec:exponentialsumapproximation}. 
Section \ref{sec:numerical_scheme} outlines the full discretization of the PDE system using implicit Runge--Kutta methods in time and finite elements in space. 
Exploiting the weak coupling of the PDE system reduces the computational cost of solving DOFPDEs to an equivalent cost of solving corresponding local-in-time PDEs. 
We apply the numerical scheme in Section \ref{sec:numerical_experiments}, recreating well-studied test cases from the literature and solving the DOFPDEs on large nontrivial 3D domains, incorporating spatial-dependency into the DO differential operator, and simulating Neo-Hookean materials with DO damping. 
We provide a conclusion in Section \ref{sec:conclusion}.

%% file: Files/KernelCompression.tex
We start by outlining kernel compression methods for FDEs and then generalize the concept to DO differential operators. 
Throughout this section, we assume that exponential sum approximations of the occurring integral kernels are given, and address how to obtain them in Section \ref{sec:exponentialsumapproximation}. 
For a recent review concerning kernel compression methods and an overview of numerical methods for FDEs, we refer the interested reader to Section 3.1 of \cite{diethelm2022trends} and \cite{karniadakis2019numerical}, respectively. 

In essence, kernel compression methods are an approximation technique for convolution operators applied to the RL integral kernel in the definition of time-fractional derivatives. Let an exponential sum approximation of $g_{1-\alpha}(t)$, $\alpha \in (0,1)$, in the form of \eqref{eq:sum_of_exponentials} be given.  
Introducing \eqref{eq:sum_of_exponentials} into the definition of the time-fractional derivative of Caputo-type \eqref{eq:def_Caputo}, we approximate the time-fractional derivative as a sum of auxiliary variables $u_k$, which we denote as fractional modes, 
\begin{equation}
    \label{eq:pre_Caputo_approximation}
    D_C^\alpha u(t) = 
    \left(g_{1-\alpha} * u^{(1)}\right)(t) 
    \approx \sum_{k=1}^m \underbrace{w_k \int_0^t e^{-\lambda_k (t-s)} u^{(1)}(s) \Diff s}_{=: u_k(t)} 
    = \sum_{k=1}^m u_k(t). 
\end{equation}
Differentiating $u_k(t)$ with respect to $t$, we find that $u_k$ is the solution of the initial value problem 
\begin{equation}
    \frac{\Diff}{\Diff t} u_k (t) = - \lambda_k  u_k(t) + w_k u^{(1)}(t), \quad  u_k(0) = 0, 
\end{equation}
and thus, obtain an approximation of $D_C^\alpha (u)$ by means of solutions of local in time differential equations. 
Numerous studies have applied kernel compression methods to solve FDEs, see, e.g., \cite{chatterjee2005statistical,diethelm2008investigation,diethelm2021new,jiang2017fast,khristenko2023solving,yuan2002numerical}. 

Now, we generalize the concept to DO differential operators. 
Naively, discretizing the integral of the fractional order and applying a kernel compression method to each time-fractional derivative is infeasible due to the large number of fractional modes. 
To achieve double precision accuracy of the kernel approximations, in the following sections, we consider up to over $1000$ quadrature points, and $60$ to $80$ exponential terms for each RL integral kernel, resulting in a total number $60,000$ to $80,000$ auxiliary variables, and thus, unreasonable memory requirements.

Similarly, we split the integral over $\alpha$ with respect to occurring integer-order derivatives, and discretize each sub-integral using a quadrature formula 
\begin{align}
    \label{eq:kc_quadrature_formula}
    \int_0^{\alphamax} \phi(\alpha) \tfC{\alpha} u(t) \Diff \alpha 
    = \sum_{i=1}^{\alphamax} \int_{i-1}^{i} \phi(\alpha) \tfC{\alpha} u(t) \Diff \alpha  
    \approx \sum_{i=1}^{\alphamax} \sum_{k=1}^{N_i} \gamma_{i,k} \phi(\alpha_{i,k}) \tfC{\alpha_{i,k}} u(t),   
\end{align}
whereby $\gamma_{i,k}$ and $\alpha_{i,k}$, $k = 1, \ldots, N_i$, are the weights and support points of an abstract open quadrature rule on the interval $(i-1, i)$. 
But instead of considering the result as multiple time-fractional derivatives, we treat the linear combinations of RL integral kernels with order $\alpha_k \in (i-1, i)$ as a single integral kernel $K_i$, 
\begin{align}
    \label{eq:kc_integral_kernel}
    \sum_{i=1}^{\alphamax} \sum_{k=1}^{N_i} \gamma_{i,k} \phi(\alpha_{i,k}) \tfC{\alpha_{i,k}} u(t) 
    = \sum_{i=1}^{\alphamax} \int_0^t \underbrace{\sum_{k=1}^{N_i} \gamma_{i,k} \phi(\alpha_{i,k}) \frac{(t-s)^{i-\alpha_{i,k}-1}}{\Gamma(i - \alpha_{i,k})}}_{=:K_i(t-s)} u^{(i)} (s) \Diff s. 
\end{align}
Note that by construction the exponents $i-\alpha_{i,k}-1 \in (-1,0)$ and the arguments of the Gamma function $i-\alpha_{i,k} \in (0,1)$ for all $k \in 1, \ldots, N_i$, $i \in 1,\ldots,\alphamax$, and thus, $K_i \in L_1[0,T]$ for any fixed $T>0$.  
We proceed as for standard time-fractional derivatives. 
Let the exponential sum approximations of the integral kernels $K_i$ be given as 
\begin{equation}
    \label{eq:kc_exponential_sum_approximations}
    K_i(t) \approx \sum_{j=1}^{m_i} w_{i,j} e^{-\lambda_{i,j}t}, \quad i = 1, \ldots, \alphamax. 
\end{equation}
Inserting \eqref{eq:kc_exponential_sum_approximations} into the approximation of the DO differential operator \eqref{eq:kc_integral_kernel}, we obtain 
\begin{align}
    \int_0^{\alphamax} \phi(\alpha) \tfC{\alpha} u(t) \Diff \alpha 
    & \approx \sum_{j=1}^{\alphamax} \sum_{i=1}^{m_j} \underbrace{w_{i,j} e^{-\lambda_{i,j}t} * u^{(i)}}_{=:u_{i,j}(t)}
    = \sum_{j=1}^{\alphamax} \sum_{i=1}^{m_j} u_{i,j}(t), 
\end{align}
whereby again the fractional modes $u_{i,j}$ are the solutions of initial value problems 
\begin{equation}
     \frac{\Diff}{\Diff t}  u_{i,j}(t) + \lambda_{i,j} u_{i,j}(t) = w_{i,j} u^{(j)}(t), \quad 
    u_{i,j}(0) = 0, \quad 
    j=1, \ldots, m_i, \quad  
    i = 1, \ldots, \alphamax. 
\end{equation}
Thus, we obtain an approximation of the DO differential operator as a sum of fractional modes given as solutions of linear, local in time differential equations. 
In contrast to the naive approach, we only require the exponential sum approximations of $\alphamax$ integral kernels. 
Since in the first step in \eqref{eq:kc_quadrature_formula}, the DO differential operator is discretized in terms of multiple time-fractional derivatives, our approach can equally be applied to multi-term FDEs.

To conclude this section, we outline the application of the kernel compression method to DOFPDEs. 
Let us consider the case 
 \begin{equation}
    \label{eq:DOFPDE}
    F\left(t, \Vec{x}, u^{(0)}(t, \Vec{x}), \ldots , u^{(n)}(t, \Vec{x}), \int_0^{\alphamax} \phi(\alpha) \tfC{\alpha} u(t, \Vec{x}) \Diff \alpha \right) = 0,  
\end{equation}
whereby the operator $F$ may depend on time $t$, the spatial variable $\Vec{x}$, and may include spatial differential operators, but no further time derivatives. 
We apply the kernel compression method to approximate the DO differential operator in terms of fractional modes  
\begin{equation}
    \label{eq:approximation_DO_operator}
    \int_0^{\alphamax} \phi(\alpha) \tfC{\alpha} u \Diff \alpha 
    \approx \sum_{i=1}^\alphamax K_i * u^{(i)}
    \approx \sum_{j=1}^{\alphamax} \sum_{i=1}^{m_j} w_{i,j} e^{-\lambda_{i,j}t} * u^{(i)}, 
\end{equation}    
and obtain the following system of PDEs 
\begin{alignat}{4}
    \label{eq:KC_constraint}
    & 0 = F \left(t, \Vec{x}, v^{(0)}, \ldots , v^{(n)}, \sum_{i=1}^{\alphamax} \sum_{j=1}^{m_i} v_{i,j}\right), \\ 
    \label{eq:KC_derivatives}
    & \frac{\Diff}{\Diff t} v^{(i-1)} 
    = v^{(i)}, \quad   
    && 
    &&& i = 1, \ldots, \alphamax, \\ 
    \label{eq:KC_modes}
    & \frac{\Diff}{\Diff t} v_{i,j} 
    = - \lambda_{i,j} v_{i,j} + w_{i,j} v^{(j)},  \quad  
    && j = 1, \ldots, m_i, \quad 
    &&& i = 1, \ldots, \alphamax.
\end{alignat}

Replacing the integral kernels with exponential sum approximations introduces a modelling error that should be balanced with the discretization errors of the numerical schemes for the temporal and spatial operators. 
We suggest a suitable numerical scheme for solving \eqref{eq:KC_constraint}-\eqref{eq:KC_modes} in Section \eqref{sec:numerical_scheme}.

%% file: Files/ExponentialSums.tex
Here, we present a constructive method to determine the weights and poles, based on the rational approximation technique \cite{khristenko2023solving} and its adaptations in \cite{beddrich2025numerical,duswald2024finite}.

As outlined in Section \ref{sec:kernelcompression}, we consider an integral kernel $K$ that arises from the discretization of the integral in \eqref{eq:kc_integral_kernel} with a quadrature formula and seek a sum of exponentials approximation $\tilde{K}$ of $K$ 
\begin{align}
    \label{eq:approximation_kernel}
    K(t) := \sum_{k=1}^N \frac{\gamma_{k} \phi(\alpha_{k})}{\Gamma(1 - \pQF{k})} t^{- \pQF{k}}, \quad  
    \tilde{K}(t) := \sum_{k=1}^m w_k e^{-\lambda_k t}, 
\end{align} 
where $\gamma_k$ and $\alpha_k$, $k=1, \ldots, N,$ denote the quadrature rule's weights and support points. 
Note that by construction, the quadrature points are contained in the interval $(0,1)$. 
Without loss of generality, we assume $\phi$ to be positive. Otherwise, we split $\phi$ into its positive and negative parts and treat these separately as two integral kernels. 
As the RL integral kernel is a completely monotone function, so is $K$, and also $\tilde{K}$ by Bernstein's theorem on completely monotone functions \cite{bernstein1929fonctions} provided that the weights $w_k$ and poles $\lambda_k$ are positive, which we numerically observe for all test cases independent of the number of exponential terms. 

\input{Files/schematic_RA}

The core idea of our approach consists of the following steps: the rational approximation of the Laplace transform of the integral kernel, a subsequent partial fractions decomposition, and obtaining a sum of exponential approximation by the inverse Laplace transform \cite{khristenko2023solving}. 
The method is based on the following standard results, which we repeat here for the reader's convenience. 
For $\alpha < 1$, the Laplace transform of the RL integral kernel can be expressed analytically, as for $\alpha < 1$, 
\begin{equation}
    \mathcal{L}\left[ t^{-\alpha} \right](t) = \Gamma(1-\alpha)s^{\alpha-1}. 
\end{equation}
Thus, the Laplace transform of $K$ is given by  
\begin{equation}
    \mathcal{L}[K](s) = \sum_{k=1}^N \omega_k \phi(\pQF{k}) s^{\alpha_k-1}, 
\end{equation}
while an exponential sum is given as the inverse Laplace transform of a partial fractions decomposition 
\begin{equation}
    \mathcal{L}^{-1}\left[ \sum_{j=1}^m \frac{w_j}{s + \lambda_j}\right](t) = \sum_{j=1}^m w_j e^{-\lambda_j t}, 
\end{equation}
whereby $w_j, \lambda_j \in \C$, for $j = 1, \ldots, m$.
A partial fractions decomposition $r(s)$ is always a rational function of type $(m-1,m)$, for a $m \in \N$, i.e., $r(s) = p(s) / q(s)$, where $p$ and $q$ are polynomials of degree $m-1$ and $m$, respectively \cite{trefethen2013approximation}. 
Accordingly, to obtain a sum of exponentials approximation, we seek a rational approximation of $\mathcal{L}[K]$ of type $(m-1,m)$. 
To compute the rational approximation, we employ the AAA algorithm \cite{nakatsukasa2018aaa}, which inherently produces a rational function of type $(m-1,m-1)$. 
Therefore, we first multiply the target function by $s$ and apply the AAA algorithm to $s\mathcal{L}[K](s)$. 
By dividing the rational approximant by $s$, we achieve a rational approximant of type $(m-1,m)$, which has a partial fraction representation \cite{trefethen2013approximation}. 
Finally, taking the inverse Laplace transform, we obtain a sum of exponentials representation. 
A sketch of the algorithmic approach is given in Figure \ref{fig:kernel_compression_schematic}.

\begin{figure}
    \centering
    \includegraphics[width=\textwidth]{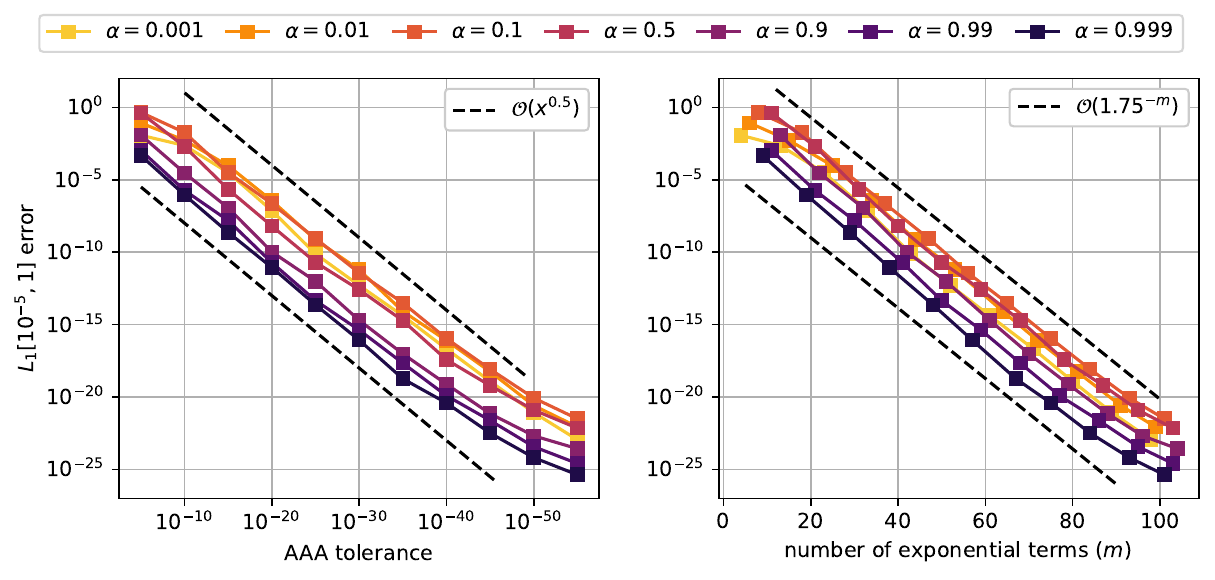}
    \caption{Error of the sum of exponential approximation of the RL integral kernels for various $\alpha$ with respect to the tolerance set for the AAA algorithm and the number of exponential terms.}
    \label{fig:CO_kernel_approximation}
\end{figure}

For standard constant order time-fractional derivatives, employing a tolerance of $10^{-13}$ for the AAA algorithm results in exponential sum approximations with up to $m = 24$ terms, whereby the $L_1([10^{-5}, 1])$-error of the kernel approximation ranges from $10^{-4}$ to $10^{-8}$ \cite{beddrich2025numerical,khristenko2023solving}, depending on the fractional order. 
Using double precision accuracy, neither reducing the tolerance of the AAA algorithm nor considering more exponential terms leads to an improvement of the $L_1([10^{-5}, 1])$-error. 
Naturally, the quality of the sum of exponentials approximation is limited by the accuracy of the underlying rational approximation. 
Given a set of discrete points and function values, the AAA algorithm results in the barycentric representation of the rational approximant with interpolation at a set of support points selected in a greedy-type fashion \cite{nakatsukasa2018aaa}. 
The quality of the rational approximation is limited by the availability and the algorithm's ability to identify suitable additional support points, whereby the addition of unsuitable support points can lead to unwanted complex poles in the rational approximation and Fréchet doublets. 
We find that, in double precision, the singular value decomposition in the AAA algorithm becomes the bottleneck of the algorithm, as no appropriate additional support points can be identified. 
Inspired by \cite{huybrechs2023aaa}, we utilize an arbitrary-precision version of the AAA algorithm \cite{driscoll2023rational} implemented in Julia \cite{bezanson2017julia}. 
Switching to arbitrary precision calculations only for the kernel approximation, we achieve double precision computational accuracy in terms of the $L_1([10^{-5}, 1])$ errors, for arbitrary RL fractional integral kernels using less than $m=80$ exponential terms.
The dependency of the $L_1([10^{-5}, 1])$-error on the tolerance given to the AAA algorithm and the number of exponential terms is shown in Figure \ref{fig:CO_kernel_approximation}. 
For the results presented, the AAA algorithm is provided a set of function values at support points
\begin{equation}
    \label{eq:AAA_support_points}
    \mathcal{X} = 
    \left\{x  \, \left| \,  x = 10^{j/25}, j = 0, 1, \ldots, 200 \right. \right\}
    \cup 
    \left\{x \, \left| \, x = 10^8 - 10^{j/25}, j = 0, 1, \ldots, 200 \right. \right\},  
\end{equation}
which is tailored to obtain a good approximation on the interval $[10^{-5},1]$. 
Heuristically, we observe that the quality of the approximation of the integral kernel is better close to $0$ if more large values are contained in $\mathcal{X}$. 
Conversely, the approximation improves for large values of $t$, the more support points close to $0$ are considered. 
This finding motivates using a logarithmic distribution of support points, which has proven superior to equidistant distributions for fractional integral kernels. 

\begin{figure}
    \centering
    \includegraphics[width=\textwidth]{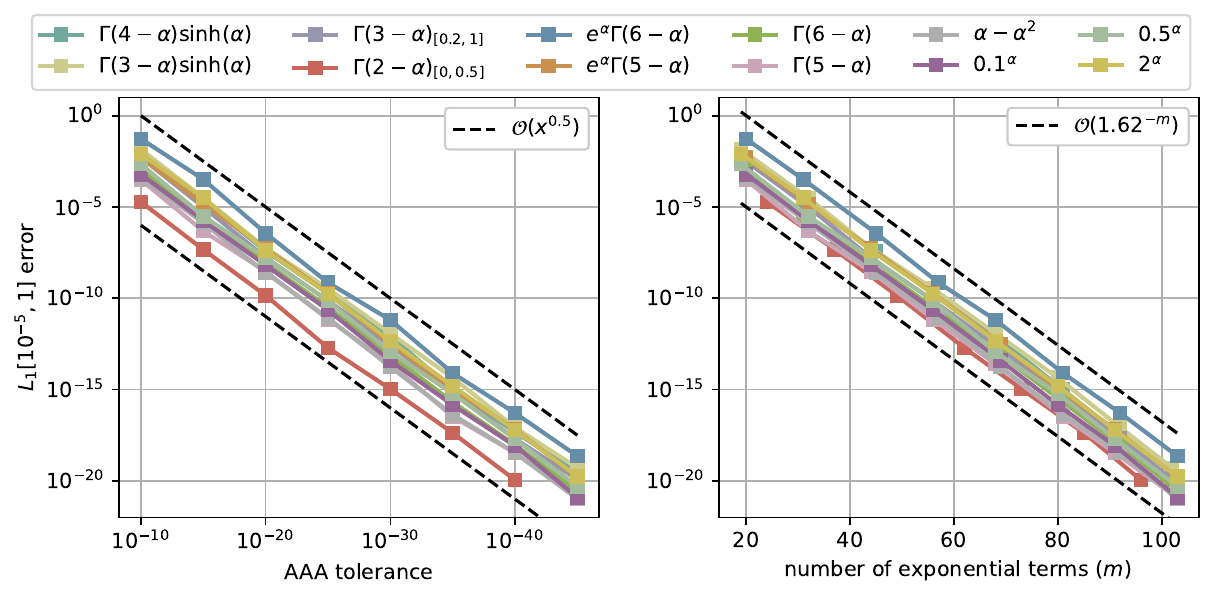}
    \caption{Error of the sum of exponential approximation integral kernels arising from several DO differential operators with respect to the tolerance set for the AAA algorithm and the number of exponential terms.}
    \label{fig:DO_kernel_approximation}
\end{figure}

In addition to the AAA-algorithm, approximating the integral kernels arising from DO differential operators requires numerical integration of the $\alpha$ integral. 
We employ an adaptive Gauss--Kronrod quadrature \cite{johnson2013quadgk} in arbitrary precision using an absolute tolerance for the quadrature of $10^{-10}$ times the tolerance set for the AAA algorithm. 
As no theoretical result ensures the convergence of the AAA algorithm, and thus, of our procedure, we demonstrate the applicability of our algorithm by approximating several integral kernels for weight functions from the literature \cite{diethelm2009numerical,katsikadelis2014numerical,yuttanan2019legendre,yuttanan2021numerical}. 
The error $\|K-\tilde{K}\|_{L_1[10^{-5},1]}$ in relation to the tolerance of the AAA algorithm and the number of exponential terms is displayed in Figure \ref{fig:DO_kernel_approximation}. 
For all considered integral kernels, we achieve an $L_1[10^{-5},1]$ error of less than $10^{-16}$ using a tolerance of $10^{-40}$ for the AAA algorithm, resulting in $m=92$ or less exponential terms. 
We attribute the higher number of exponential terms for the DO integral kernels than for the RL integral kernels to the scaling introduced by $\phi(\alpha)$.

Note that the additional computational cost of calculating the exponential sum approximation in arbitrary precision is part of the offline phase of a simulation and does not influence the cost of computing a single time step. While the rational approximation and partial fractions decomposition are performed with arbitrary precision, we store and reuse the weights and poles with double precision. 

As the AAA algorithm can be applied to an arbitrary set of support points and corresponding function values, the presented approach is independent of the specific form of $K$. 
Indeed, it is not even necessary to obtain an analytical expression of the Laplace transform of the integral kernel, but it can be evaluated numerically. 
Therefore, the method is well-suited for generalizations to more complex integral kernels, and can also be applied to tempered fractional integral kernels \cite{sabzikar2015tempered}.

%% file: Files/schematic_RA.tex
\tikzset{
    modernbox/.style={
        draw,
        rectangle,
        thick,
        rounded corners=8pt,
        fill=blue!15,
        inner sep=10pt,
        minimum height=1.2cm,
        drop shadow
    },
    arrowstyle/.style={
        ->,
        thick,
        >= Triangle
    },
    dashedarrow/.style={
        ->,
        thick,
        dashed,
        >=Triangle
    }
}
\usetikzlibrary{shadows, arrows.meta}
\begin{figure}
    \centering
    \resizebox{10cm}{5cm}{\usetikzlibrary{shadows, arrows.meta}
    \begin{tikzpicture}
    \node[modernbox] (sum1) at (0,6) 
        {\Large$\displaystyle \sum_{k = 1}^N \JB{\gamma}_{k} \phi(\pQF{k}) s^{\alpha_k}$};
    \node[modernbox] (sum2) at (0,3) 
        {\Large$\displaystyle \sum_{k = 1}^N \JB{\gamma}_{k} \phi(\pQF{k}) s^{\alpha_k-1}$};
    \node[modernbox] (sum3) at (0,0) 
        {\Large$\displaystyle \sum_{k = 1}^N \frac{\JB{\gamma}_k \phi(\pQF{k})}{\Gamma(1 - \pQF{k})} t^{- \pQF{k}}$};
    
    \node[modernbox] (frac1) at (5,6) 
        {\Large$\displaystyle \frac{p(s)}{q(s)}$};
    \node[modernbox] (frac2) at (5,3) 
        {\Large$\displaystyle \frac{p(s)}{s q(s)}$};

    \node[modernbox] (frac3) at (10,3) 
        {\Large$\displaystyle \sum_{j=1}^m \frac{w_j}{s+\lambda_j}$};
    \node[modernbox] (frac4) at (10,0) 
        {\Large$\displaystyle \sum_{j=1}^m w_j e^{-\lambda_j t}$};

    \draw[arrowstyle] (sum2) -- (sum1) node[midway, right, xshift=3pt] {\LARGE$\cdot s$};
    \draw[arrowstyle] (sum3) -- (sum2) node[midway, right, xshift=3pt] {\LARGE$\mathcal{L} \{\cdot\}(s)$};
    \draw[arrowstyle] (sum1) -- (frac1) node[midway, above, yshift=4pt] {\LARGE AAA};
    \draw[arrowstyle] (frac1) -- (frac2) node[midway, right, xshift=4pt] {\LARGE$\cdot / s$};
    \draw[arrowstyle] (frac2) -- (frac3) node[midway, above, yshift=4pt] {\LARGE PFD};
    \draw[arrowstyle] (frac3) -- (frac4) node[midway, right, xshift=3pt] {\LARGE$\mathcal{L}^{-1} \{\cdot \}(t)$};
    \draw[dashedarrow] (sum3) -- (frac4) node[midway, above, yshift=4pt] {\LARGE$\approx$};
    \end{tikzpicture}
    }
    \caption{
    Schematic sketch of the kernel compression approach: 
    (1) Calculation of the analytical Laplace transform, 
    (2) Multiplication with $s$, 
    (3) Rational approximation using the AAA algorithm, 
    (4) Division by $s$, 
    (5) Calculation of the partial fractions decomposition, 
    (6) Calculation of the analytical inverse Laplace transform.}
    \label{fig:kernel_compression_schematic}
\end{figure}

%% file: Files/NumericalMethod.tex
This section addresses the spatial and temporal discretization of the PDE systems of \eqref{eq:KC_constraint}-\eqref{eq:KC_modes}. 
Application of the kernel compression method to DOFPDEs results in systems of PDEs, whereby the number of fractional modes depends on the accuracy of the approximation of the integral kernels. 
Exploiting that the time derivative of each fractional mode is given by a linear combination of itself and a time derivative of the solution, we employ static condensation to reduce the PDE system to an implicit problem. 
We emphasize that the size of the implicit problem is independent of the number of exponential terms. 
As a result, the solution of the DOFPDE can be calculated at an equivalent computational cost to solving a corresponding local-in-time PDE. 

\subsection{Temporal discretization}

We apply Runge--Kutta methods to the PDE system, whereby the stiffness of the fractional mode equations restricts us to implicit methods. 
Further, the wide range of the poles of the exponential sum approximations \cite{khristenko2023solving} results in both slow and fast decaying components, necessitating the use of an L-stable method, see, e.g., Section 6.1 in \cite{khristenko2023solving}. 
Therefore, we consider the RadauIIA Runge--Kutta schemes, as well as suitable Diagonally Implicit Runge--Kutta (DIRK) schemes and refer to the textbooks \cite{brenan1995numerical,hairer1989runge,hairer1993solving} for detailed introductions. 
For readability, we omit the dependency on the spatial variable $\Vec{x}$, $j = 1, \ldots, m_i$, and $i=1,\ldots, \alphamax$ throughout this section.

To compute the discrete values $v^{(i),n+1}$ and $v_{i,j}^{n+1}$ at $t^{n+1}$, we apply an $s$-stage Runge--Kutta method, identified by the coefficients $\tensor{A} \in \R^{s\times s}$, $\Vec{b}, \Vec{c} \in \R^s$, to \eqref{eq:KC_constraint}-\eqref{eq:KC_modes}. 
Given discrete values at $t^n$, the vectors of Runge--Kutta steps, denoted by $\Vec{k^{(i)}}$ and $\Vec{k^{i,j}}$ for the variables $v^{(i)}$ and $v_{i,j}$, are implicitly given by 
\begin{alignat}{4}
    \label{eq:RK_constraint}
    0 = F \Bigg(\Vec{t^{n+c}}, 
    v^{(0),n} \Vec{1} & + h \tensor{A} \Vec{k^{(0)}}, \ldots , v^{(\alphamax),n} \Vec{1} + h \tensor{A} \Vec{k^{(\alphamax)}}, 
    \sum_{i=1}^{\alphamax} \sum_{j=1}^{m_i} \left(v_{i,j}^{n} \Vec{1} + h \tensor{A} \Vec{k^{i,j}} \right)\Bigg), \\ 
    \label{eq:RK_derivatives}
    \Vec{k^{(i-1)}} 
    & = v^{(i),n} \Vec{1} + h \tensor{A} \Vec{k^{(i)}}, \\ 
    \label{eq:RK_modes}
    \Vec{k^{i,j}} 
    & = - \lambda_{i,j} \left(v_{i,j}^{n} \Vec{1} + h \tensor{A} \Vec{k^{i,j}} \right)+ w_{i,j} \Vec{k^{(i-1)}},
\end{alignat}
whereby $\Vec{t^{n+c}} = (t^n + c_1 h, \ldots, t^n + c_s h)\Transpose$ and, for vector inputs, the operator $F$ is applied componentwise.
We recall that for RadauIIA schemes and fully implicit DIRK schemes, $\tensor{A}$ has only positive eigenvalues. 
In combination with $\lambda_{i,j} \geq 0$, it is obvious, that the matrix $\tensor{B^{i,j}} = \left(\tensor{Id} + \lambda_{i,j} h \tensor{A}\right)^{-1}$ is well-defined. 
Thus, we can solve \eqref{eq:RK_modes} for $\Vec{k^{i,j}}$ and obtain 
\begin{align}
    \label{eq:RK_modes_solved}
    \Vec{k^{i,j}}  
    & = \tensor{B^{i,j}} \left(- \lambda_{i,j} v_{i,j}^n \Vec{1}
    + w_{i,j} \Vec{k^{(i-1)}} \right). 
\end{align}
Repeatedly applying \eqref{eq:RK_derivatives}, we can express $\Vec{k^{(i)}}$ in terms of $\Vec{k^{(\alphamax)}}$ 
\begin{align}
    \label{eq:RK_derivatives_solved}
    \Vec{k^{(i)}} = \sum_{l=1}^{\alphamax-i} \left( h \tensor{A}\right)^{l-1} \Vec{1} v^{(i+l),n} + \left( h \tensor{A}\right)^{\alphamax-i} \Vec{k^{(\alphamax)}}. 
\end{align}
Inserting \eqref{eq:RK_modes_solved} and \eqref{eq:RK_derivatives_solved} into \eqref{eq:RK_constraint}, we obtain an implicit condensed system soley dependent on $\Vec{k^{(\alphamax)}}$
\begin{equation}
    \label{eq:RK_condensed_system}
    0 = F \Bigg(\Vec{t^{n+c}}, 
    \Vec{\mathrm{v}^{(0)}}\left(\Vec{k^{(\alphamax)}}\right), \ldots , 
    \Vec{\mathrm{v}^{(\alphamax)}}\left(\Vec{k^{(\alphamax)}}\right), 
    \sum_{i=1}^{\alphamax} \sum_{j=1}^{m_i} \Vec{\mathrm{v}_{i,j}}\left(\Vec{k^{(\alphamax)}}\right) \Bigg), 
\end{equation}
whereby 
\begin{alignat}{4}
    \label{eq:RK_condensed_system_notation}
    \Vec{\mathrm{v}^{(i)}}\left(\Vec{k^{(\alphamax)}}\right) & = v^{(i),n} \Vec{1} + \sum_{l=1}^{\alphamax-i} v^{(i+l),n} \left( h \tensor{A}\right)^{l} \Vec{1}  + \left( h \tensor{A}\right)^{\alphamax-i+1} \Vec{k^{(\alphamax)}}, \\ 
    \Vec{\mathrm{v}_{i,j}}\left(\Vec{k^{(\alphamax)}}\right) & = v_{i,j}^{n} \Vec{1} 
    - h \lambda_{i,j} v_{i,j}^n \tensor{A} \tensor{B^{i,j}} \Vec{1} \\ \nonumber
    & + h w_{i,j} \tensor{A} \tensor{B^{i,j}} \left( \sum_{l=1}^{\alphamax-i+1}v^{(i-1+l),n}  \left( h \tensor{A}\right)^{l-1} \Vec{1} + \left( h \tensor{A}\right)^{\alphamax-i+1} \Vec{k^{(\alphamax)}} \right). 
\end{alignat}
The other vectors of Runge--Kutta steps $\Vec{k^{(i)}}$ and $\Vec{k^{i,j}}$ can be explicitly determined by a linear update using \eqref{eq:RK_modes_solved} and \eqref{eq:RK_derivatives_solved}, and we recover the discrete values of the solution, its time derivatives, and the fractional modes by  
\begin{equation}
    v^{(i),n+1} = v^{(i),n} + h \Vec{b}\Transpose \Vec{k^{(i)}}, \quad 
    v_{i,j}^{n+1} = v_{i,j}^{n} + h \Vec{b}\Transpose \Vec{k^{i,j}}. 
\end{equation}
It is essential to point out that the condensed system \eqref{eq:RK_condensed_system} results in an $s \times s$ block structure, the same structure as for applying an $s$-stage Runge-Kutta method to solve an integer-order PDE. 
Therefore, the number of exponential terms of the kernel approximation only influences the number of linear updates of the fractional modes, not the size of the implicit problem \eqref{eq:RK_condensed_system}. 

\subsection{Spatial discretization}

The PDE systems are discretized in space using second-order, continuous finite elements. 
Given a possibly unstructured mesh, we use Gauss--Lobatto basis functions to define the finite element spaces and recover the coefficient vectors representing the approximations of the solutions in terms of the polynomial basis functions. 
We must store $\alphamax + \sum_{i=1}^\alphamax m_i$ coefficient vectors for the solution, its time derivatives, and the fractional modes. 
As outlined in the previous section, after application of an $s$-stage Runge--Kutta scheme, the PDE system simplifies to a condensed system of $s$ coupled PDEs that has to be solved implicitly. 

At the same time, the coefficient vectors of the fractional modes are efficiently computed by a linear update without the need to solve a PDE, inverting a mass-matrix, or even a matrix-vector multiplication.
Consequently, no significant computational cost arises for numerous terms in the exponential sum approximation. 
Despite the number of fractional modes depending on the accuracy of the kernel approximation, it is therefore not essential to balance the modelling error introduced by the kernel approximation with the error of the temporal and spatial discretizations, but it can always be chosen sufficiently small.
Ultimately, our numerical scheme for DOFPDE carries no significant additional cost in comparison to solving corresponding, local-in-time PDEs.   

The stability and convergence of fully implicit Runge--Kutta methods for PDEs entail computational implications in the form of highly coupled, poorly conditioned linear systems \cite{mardal2007order}. 
For the numerical examples in Section \ref{sec:numerical_experiments}, a general algebraic multigrid preconditioner proves to be sufficient, and we suggest the use of DIRK schemes as a possible alternative to avoid poorly conditioned systems. While the following numerical examples are realized using the MFEM library \cite{mfem}, our approach can be easily implemented in any finite element library. 

%% file: Files/NumericalResults.tex
We apply the kernel compression method and numerical scheme introduced in Sections \ref{sec:kernelcompression}-\ref{sec:numerical_scheme} to various DOFDEs and DOPDEs. 
First, we recreate well-studied test cases from the literature before solving linear and nonlinear DOFPDEs. 

\subsection{Numerical example 1: Influence of the kernel approximation}
\label{subsec:num_ex_kernel_approximation}
\begin{figure}
    \centering
    \includegraphics[width=\textwidth]{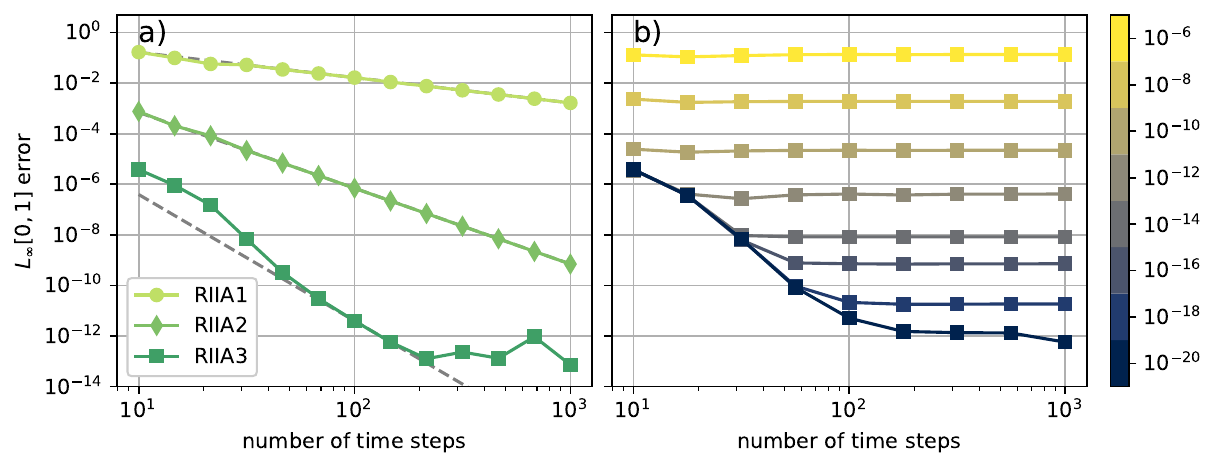}
    \caption{Numerical example 1: (a) Convergence study with respect to the number of time steps. The system of ODEs is solved using the $s$-stage Radau IIA scheme (RIIA$s$), $s=1,2,3$, with expected convergence order $2s-1$. The kernel is approximated using a tolerance of $10^{-40}$ for the AAA algorithm. (b) Influence of the rational approximation quality quantified in terms of the tolerance for the AAA algorithm, ranging from $10^{-6}$ to $10^{-20}$. The system of ODEs is solved using the RIIA3 scheme.}
    \label{fig:Ex1_combined}
\end{figure}

We consider a well-studied DOFDE from the literature \cite{diethelm2009numerical}. 
Let $\phi:[0,2]\rightarrow \R$, $\phi(\alpha) = e^{-\alpha}\Gamma(6-\alpha)$, and initial conditions $u(0) = u^{(1)}(0) = 0$ be given. Then $u(t) = t^5$ is the analytical solution of 
\begin{equation}
    \int_0^2 \phi(\alpha) D^\alpha_Cu(t)\Diff \alpha 
    = 120 \frac{t^5 - t^3/e^2}{1+\log(t)}. 
\end{equation}
In contrast to the original setting, where the right-hand side is solely dependent on $t$, through the addition of $g(u(t))$ and subtraction of $g(t^5)$, we obtain a solution-dependent right-hand side, while preserving the analytical solution. 
Choosing $g(u) = \cos(100u / (1 + \sqrt{u}))$, the problem reads as  
\begin{equation}
    \int_0^2 \phi(\alpha) D^\alpha_Cu(t)\Diff \alpha = 120 \frac{t^5 - t^3/e^2}{1+\log(t)} 
    + \cos\left(\frac{100u}{1+\sqrt{u}}\right)
    - \cos\left(\frac{100t^5}{1+\sqrt{t^5}}\right). 
\end{equation}
The convergence of the $L_\infty[0,1]$-error is displayed in Figure \ref{fig:Ex1_combined}a. 
The system of ordinary differential equations is solved using an $s$-stage RadauIIA (RIIA$s$) scheme, $s=1,2,3$. 
After initial pre-asymptotic behavior, the convergence rates match the expected order of $2s-1$. The occurring nonlinear system is solved up to an absolute error of $10^{-12}$, resulting in a plateau for the RIIA3 scheme and time steps smaller than $1/200$. 
For the results presented in Figure \ref{fig:Ex1_combined}a, the tolerance of the AAA algorithm is chosen as $10^{-40}$, resulting in $92$ and $91$ exponential terms for the two integral kernels, respectively.  
The error introduced by the exponential sum approximations of the integral kernels is dominated by the tolerance of the nonlinear solver, and thus, it is not visible in the figure. 
Therefore, we further study the influence of the kernel approximation in Figure \ref{fig:Ex1_combined}b, where we consider a range of exponential sum approximations of the integral kernels with tolerances for the AAA algorithm from $10^{-6}$ to $10^{-20}$, corresponding to $10$ to $44$ exponential terms per integral kernel, such that other error terms do not dominate the approximation error. 
The lower the chosen tolerance, the smaller the level at which the $L_\infty[0,1]$-error plateaus. 

\subsection{Numerical example 2: Non-zero initial conditions}
\label{subsec:non_zero_initial_conditions}

For the second example, we consider non-zero initial conditions \cite{diethelm2009numerical}, which poses the challenge that independently of the regularity of the right-hand side, the solution is possibly non-smooth at $t=0$ \cite{diethelm2010analysis}, and thus, numerical methods do not achieve their optimal convergence order on uniform grids.  
To overcome this obstacle, we employ a graded mesh for the temporal discretization 
\begin{equation}
    \label{eq:graded_mesh}
    t_n = \left(\frac{n}{N}\right)^\gamma T, 
\end{equation}
where $N$ denotes the number of time steps and $\gamma$ the grading factor, see, e.g., \cite{stynes2017error}.
Let $\phi:[0,2]\rightarrow \R$, $\phi(\alpha) = \Gamma(4-\alpha)$, and $u(0) = 4$, $u^{(1)}(0) = 2$. 
Then $u(t) = t^3+2t+4$ is the analytical solution of the DOFDE
\begin{equation}
    \int_0^2 \phi(\alpha) D^\alpha_Cu(t)\Diff \alpha 
    = \frac{6t^3 + 6t -4 }{\log t} + \frac{6 - 10t}{(\log t)^2} + \frac{4t-4}{(\log t)^3}. 
\end{equation}
Note that the right-hand side of the problem, and thus, the solution of the multi-term FDE, which we obtain after applying the quadrature rule to the DOFDE, and the solution of the mode system, are not continuously differentiable at $t=0$. 

\begin{figure}
    \centering
    \includegraphics[width=\textwidth]{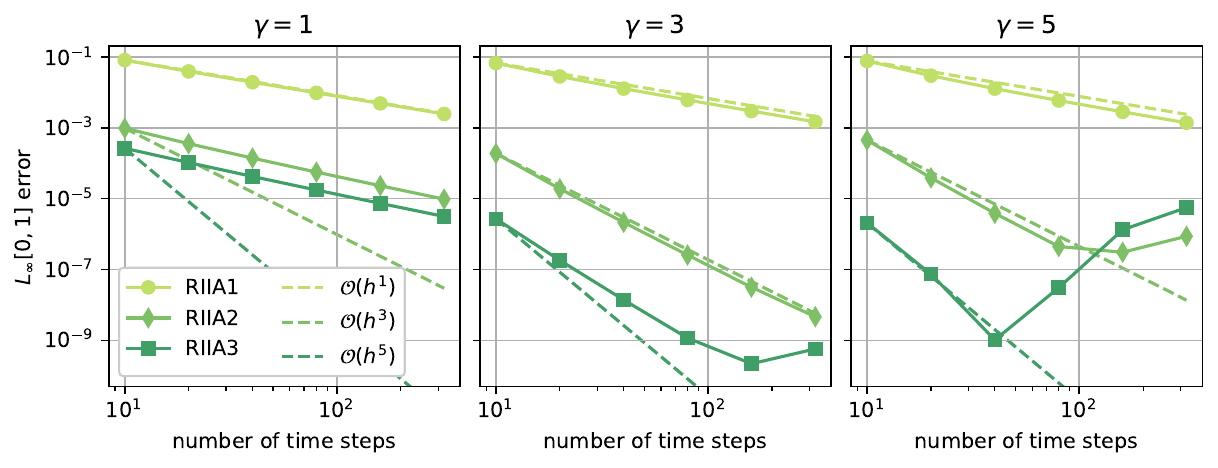}
    \caption{Numerical example 2: (a) Convergence study with respect to the number of time steps for different graded meshes with $\gamma=1,3,5$. The ODE system is solved using the $s$-stage Radau IIA scheme (RIIA$s$), $s=1,2,3$, with expected convergence of order $2s-1$ for a sufficiently smooth solution.}
    \label{fig:example2}
\end{figure}

For the kernel compression method, we set the tolerance of the AAA algorithm to $10^{-45}$, resulting in $m=102$ exponential terms for all integral kernels. 
We employ the $s$-stage Runge--Kutta method, $s=1,2,3$, using different mesh gradings $\gamma = 1,3,5$ with up to $N=320$ time steps.  
The $L_\infty[0,1]$-error is displayed in Figure \ref{fig:example2}. 
For $\gamma=1$, we observe first-order convergence of all Runge--Kutta schemes, which is expected due to the lack of regularity of the approximate system. 
The results show third-order convergence for $s=2,3$ choosing $\gamma = 3$, and fifth-order convergence for $s=3$ and $\gamma=5$. 
Due to rounding errors, for the first, very small time steps, we observe that the error plateaus significantly earlier for high values of $\gamma$.
The results indicate that the grading factor of the time steps should be chosen as the convergence order of the temporal discretization to achieve the optimal error decay while reducing spurious effects. 

\subsection{Numerical example 3: Distributed-order diffusion-wave equation}

We consider the DO diffusion-wave equation, one of the simplest and most studied PDEs in the context of DO operators. 
The system reads as  
\begin{equation}
    \label{eq:DO_diffusion_wave}
    \int_0^2 \phi(\alpha) \tfC{\alpha} u(t,\Vec{x}) \Diff \alpha = \epsilon \Laplace u(t,\Vec{x}) + f(t,\Vec{x}),  
\end{equation}
subject to initial conditions for the solution and its first derivative 
\begin{alignat}{3}
u(0,\Vec{x}) 
& = u_0(\Vec{x}), \quad 
u^{(1)}(0,\Vec{x}) 
& = v_0(\Vec{x}), \quad 
&& \Vec{x} \in \Omega, 
\end{alignat}
as well as Dirichlet and Neumann boundary conditions  
\begin{alignat}{3}
u(t,\Vec{x}) 
& = u_D(t,\Vec{x}), \quad 
&& (t,\Vec{x}) \in (0,T) \times \Gamma_D, \\ 
\nabla u(t,\Vec{x}) \cdot \Vec{n} 
& = u_N (t,\Vec{x}), \quad 
&& (t,\Vec{x}) \in (0,T) \times \Gamma_N. 
\end{alignat}

We apply the kernel compression method to the DO operator, as outlined in Section \ref{sec:kernelcompression}. 
Given a suitable sum of exponential approximations, we obtain the PDE system 
\begin{alignat}{3}
    \sum_{i=1}^2 \sum_{j=1}^{m_i} v_{i,j}(t,\Vec{x}) & = \epsilon \Laplace v(t,\Vec{x}) + f(t,\Vec{x}),  \\ 
    \frac{\mathrm{d}}{\mathrm{d}t} v_{i,j}(t,\Vec{x}) & = -\lambda_{i,j} v_{i,j}(t,\Vec{x}) + w_{i,j} v (t,\Vec{x}), \quad 
    && j = 1, \ldots, m_i, \quad 
    i = 1,2, 
\end{alignat}
subject to the same initial and boundary conditions 
\begin{alignat}{3}
v(0,\Vec{x}) 
 = u_0(\Vec{x}), \quad 
v^{(1)}(0,\Vec{x})
& = v_0(\Vec{x}), \quad 
&& \Vec{x} \in \Omega, \\ 
v(t,\Vec{x}) 
& = u_D(t,\Vec{x}), \quad 
&& (t,\Vec{x}) \in (0,T) \times \Gamma_D, \\ 
\nabla v(t,\Vec{x}) \cdot \Vec{n} 
& = u_N (t,\Vec{x}) \quad 
&& (t,\Vec{x}) \in (0,T) \times \Gamma_N. 
\end{alignat}

\subsubsection{Convergence study}

\renewcommand{\arraystretch}{0.8}
\begin{table}[h!]
\centering
\begin{tabular}{cc|cc|cc|cc}
\toprule
& & \multicolumn{2}{c|}{Implicit Euler} & \multicolumn{2}{c|}{DIRK ($s=2$)} & \multicolumn{2}{c}{Radau IIA ($s=2$)} \\  
\midrule
$\Delta x$ & $h$ & error & rate & error & rate & error & rate \\
\midrule
1/8   & 1/8 & 1.85e-2 & -- & 1.95e-2 & -- & 1.64e-2 & -- \\
1/16  & 1/16 & 1.01e-2 & 0.87 & 2.72e-3 & 2.85 & 1.98e-3 & 3.06 \\
1/32  & 1/32 & 5.37e-3 & 0.91 & 4.60e-4 & 2.56 & 2.46e-4 & 3.00 \\
1/64  & 1/64 & 2.74e-3 & 0.97 & 9.70e-5 & 2.24 & 3.07e-5 &  2.99 \\
1/128  & 1/128 & 1.38e-3 & 0.99 & 2.30e-5 & 2.07 & 3.85e-6 &  2.99 \\
1/256  & 1/256 & 6.91e-4 & 1.00 & 5.68e-6 & 2.01 & 4.81e-7 &  2.99 \\
\bottomrule
\end{tabular}
\caption{Combined convergence study of the $L^\infty([0,T]; L^2(\Omega))$-error with respect to spatial and temporal discretization.}
\label{tab:convergence}
\end{table}

Let $\phi(\alpha) = e^{-\alpha}\Gamma(6-\alpha)$, $T=1$, $\Omega = (0,1)^2$, $\epsilon=1$, and the right-hand side is given by 
\begin{equation}
    f(t,\Vec{x}) = \sin(4 \pi x_1) \sin(4 \pi x_2) 
    \left( \frac{120(t^5 - t^3e^{-2})}{1 + \log(t)} + 32 \pi^2 t^5\right). 
\end{equation}
Then $u(t,\Vec{x}) = \sin(4 \pi x_1) \sin(4 \pi x_2) t^5$ is the analytical solution of \eqref{eq:DO_diffusion_wave}, whereby the system is closed by homogeneous Dirichlet boundary conditions. 
For the rational approximation, we set the tolerance of the AAA algorithm as $10^{-20}$ for all integral kernels and choose the tolerance for the numerical quadrature and the support points as outlined in Section \ref{sec:exponentialsumapproximation}, yielding in $m_1 = 45$, $m_2 = 44$ exponential terms for the two integral kernels. 
As we employ second-order finite elements for the spatial discretization, we expect a cubic decay rate of the spatial $L^2(\Omega)$ error. 
We discretize the system in time using the implicit Euler method, an $L$-stable $2$-stage DIRK scheme, and the $2$-stage Radau IIA scheme with first, second, and third order convergence, respectively. 
The error $\|u_h -u\|_{L^\infty([0,T]; L^2(\Omega))}$ is displayed in Table \ref{tab:convergence}. 
Asymptotically, the decay matches the theoretical expectations. At the same time, we observe a better-than-expected error decay for the DIRK scheme at coarse discretizations, where the spatial error component still dominates the temporal error.

\subsubsection{A non-trivial 3D domain}

\begin{figure}
    \centering
    \begin{subfigure}[b]{0.49\textwidth}
         \centering
         \includegraphics[width=\textwidth, trim = {15cm 17cm 16cm 26cm}, clip]{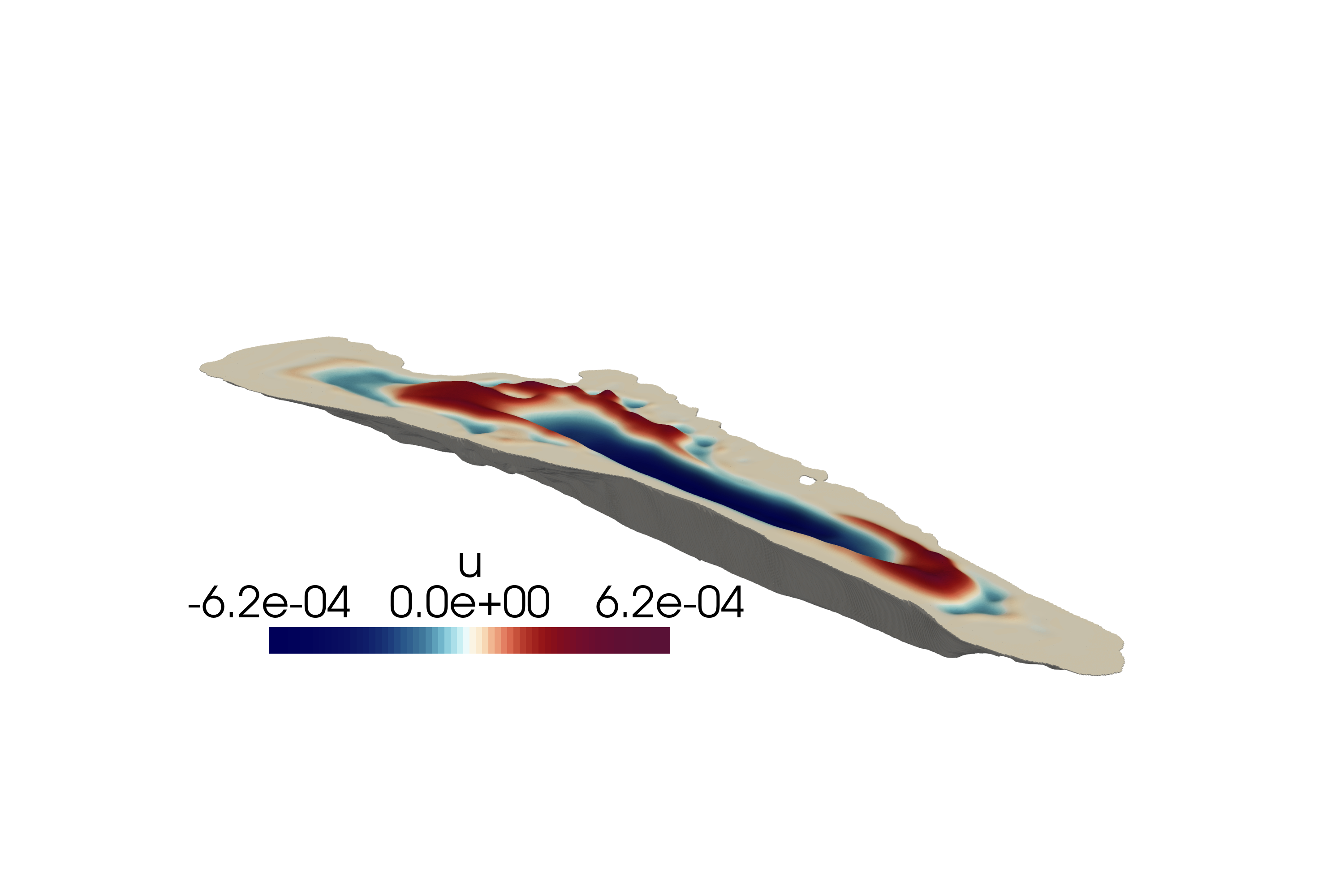}
         \caption{$\phi_{2, 0.1}(\alpha)$, $t=0.5$}
    \end{subfigure}
    \hfill
    \begin{subfigure}[b]{0.49\textwidth}
         \centering
         \includegraphics[width=\textwidth, trim = {15cm 17cm 16cm 26cm}, clip]{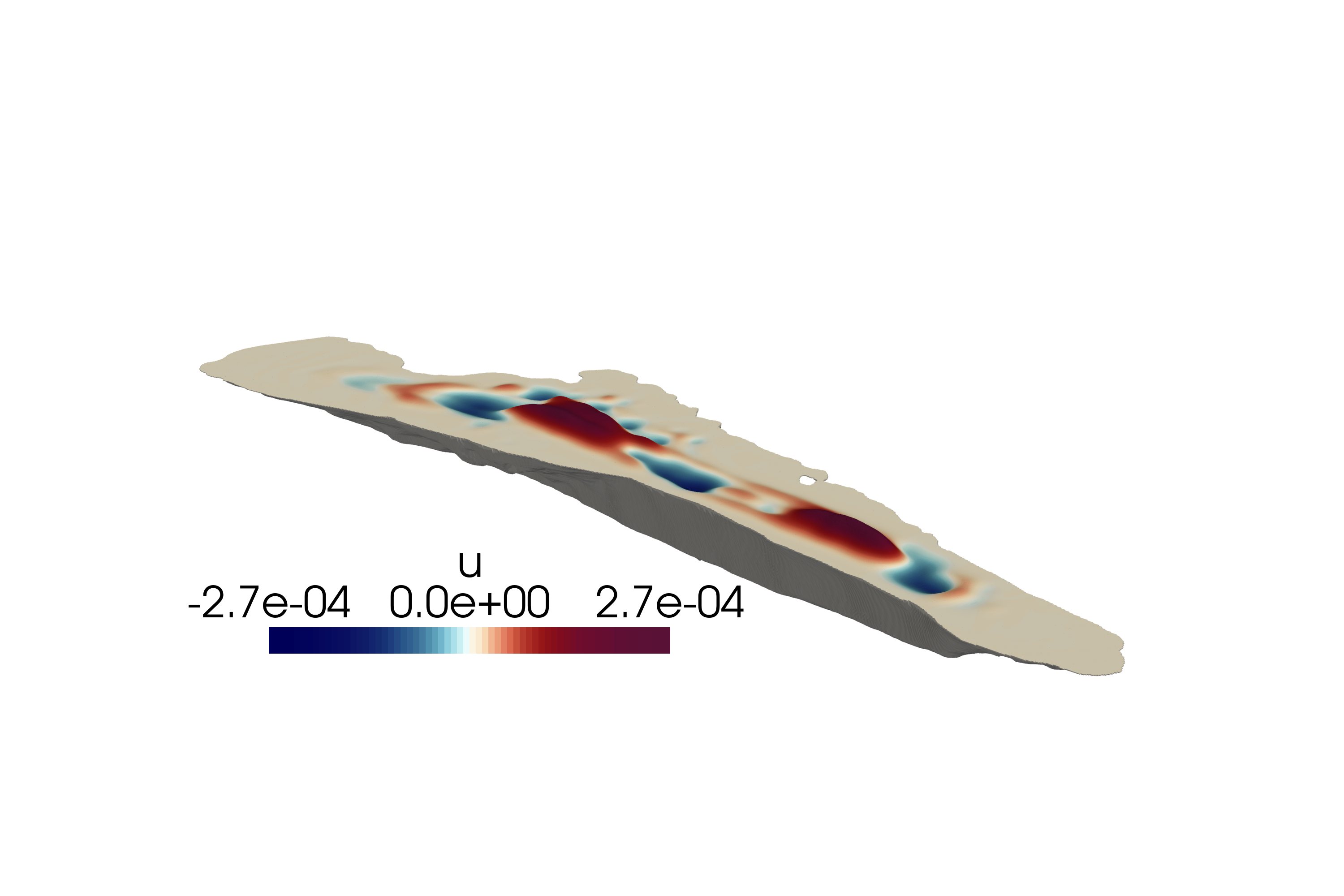}
         \caption{$\phi_{2, 0.1}(\alpha)$, $t=1$}
    \end{subfigure}
    \begin{subfigure}[b]{0.49\textwidth}
         \centering
         \includegraphics[width=\textwidth, trim = {15cm 17cm 16cm 26cm}, clip]{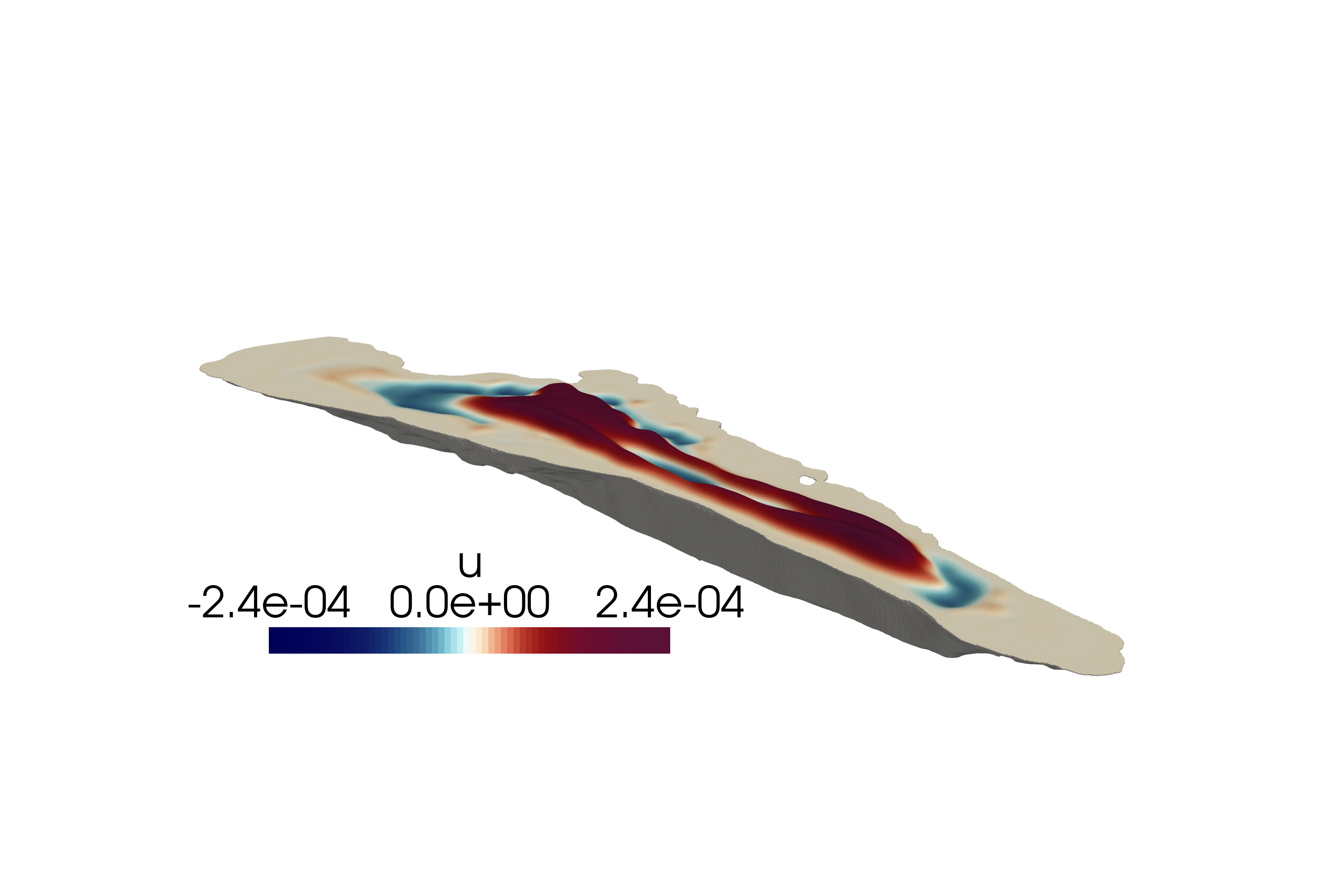}
         \caption{$\phi_{2, 0.5}(\alpha)$, $t=0.5$}
    \end{subfigure}
    \hfill 
    \begin{subfigure}[b]{0.49\textwidth}
         \centering
         \includegraphics[width=\textwidth, trim = {15cm 17cm 16cm 26cm}, clip]{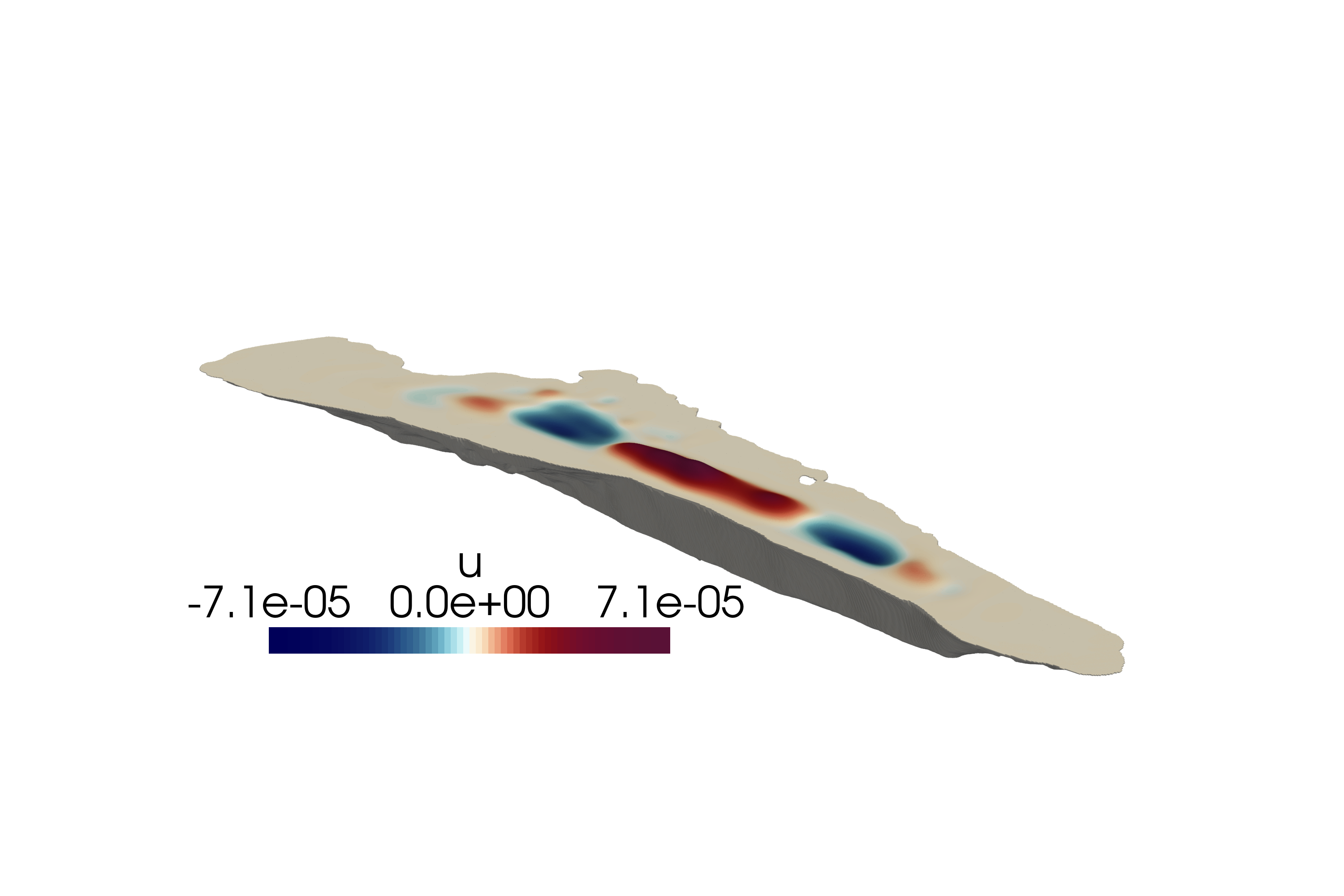}
         \caption{$\phi_{2, 0.5}(\alpha)$, $t=1$}
    \end{subfigure}
    \caption{Solution of the DO diffusion-wave equation for Lake Starnberg at $t=0.5,1$. For better visibility of the depth profile, the $z$-axis of the domain is rescaled by a factor of 10.}
    \label{fig:Starnbergersee}
\end{figure}

To demonstrate our method's robustness, we solve the DO diffusion-wave equation on a 3D geometry based on Lake Starnberg in southern Germany, whereby the mesh comprises $40.719.158$ dofs for the second-order finite element discretization.
We consider homogeneous Neumann boundary conditions for the water surface and homogeneous Dirichlet boundary conditions otherwise, and set $T=2$, $\epsilon=1000$, and $u_0 \equiv v_0 \equiv 0$. 
The external force 
\begin{equation}
    f(t, \Vec{x}) = \begin{cases}
        100 \sin(20 \pi t) * \exp\left(\frac{-1}{10t(1-10t)}\right) & t\in (0,0.1), \\ 
        0 & \text{otherwise}, 
    \end{cases}
\end{equation}
acts in the form of one oscillation uniformly on the entire domain. 
We compare the solution for two different weight functions $\phi_{2,0.1}$ and $\phi_{2,0.5}$, given as 
\begin{align}
    \phi_{2, r}(\alpha) = c_r \exp\left(\frac{1}{(\alpha-2)^2 - r^2}\right) \chi_{[2-r,2]}(\alpha), 
\end{align}
where $c_r$ denotes an $r$-dependent normalizing constant. 
While for $\phi_{2,0.1}$ the DO diffusion-wave equation closely resembles the wave equation, $\phi_{2,0.5}$ also sets a positive weight for $\alpha \in (1.5,1.9)$, strengthening the diffusive character of the DO diffusion-wave equation. The integral kernels are approximated using $m=20$ exponential terms, and the size of the timesteps is set to $h=0.01$. 

Figure \ref{fig:Starnbergersee} displays the solutions for both weight functions at $t=0.5$ and $t=1$. 
We observe the influence of the domain geometry in the form of reflections in the wave profile, as well as qualitative differences between the two scenarios.
The magnitude of the solution for $\phi_{2,0.5}$ is reduced compared to $\phi_{2,0.1}$, which is expected due to the increased diffusive behavior of the respective DO differential operator. 

\subsubsection{Space-dependent DO differential operators}

\begin{figure}
    \centering
    \begin{subfigure}[b]{\textwidth}
         \centering
         \includegraphics[width=\textwidth]{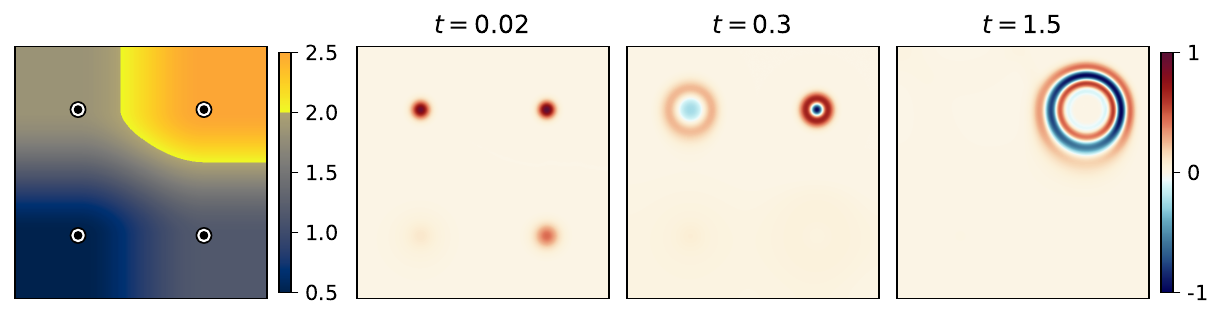}
         \caption{Geometric scenario}
         \label{fig:Geometric}
    \end{subfigure}
    \begin{subfigure}[b]{\textwidth}
         \centering
         \includegraphics[width=\textwidth]{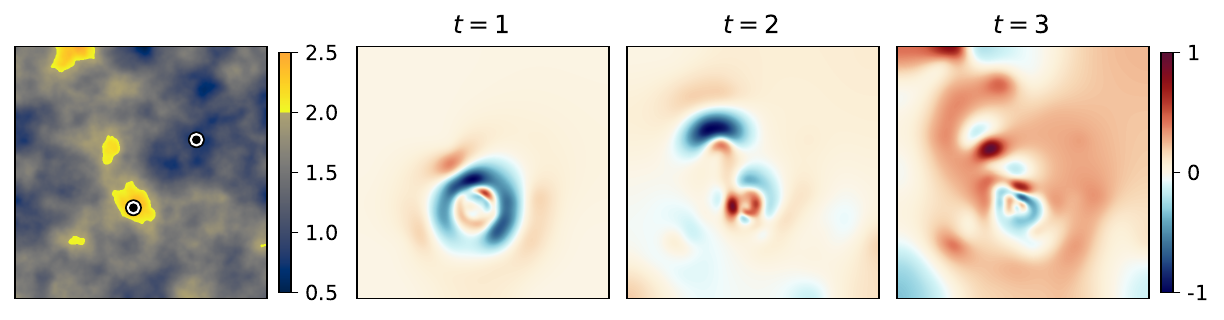}
         \caption{Random field scenario}
         \label{fig:RandomField}
    \end{subfigure}
    \begin{subfigure}[b]{\textwidth}
         \centering
         \includegraphics[width=\textwidth]{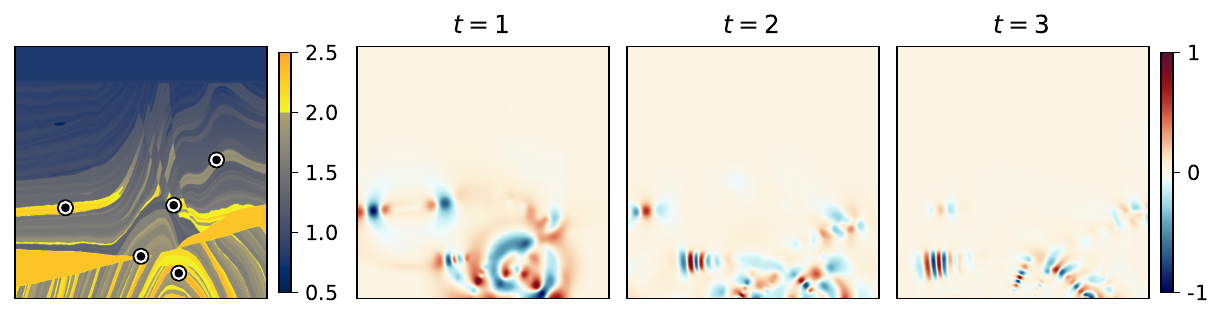}
         \caption{Marmousi scenario}
         \label{fig:Marmousi}
    \end{subfigure}    
    \caption{Spatial dependency of $\eta(\Vec{x})$ (left).
    Areas where $\eta(\Vec{x})>2$ are highlighted, and the positions of the point sources of the initial conditions are visualized as black and white circles. 
    Normalized solution of the DO diffusion-wave equation with a spatially dependent weight function for different $t$ (middle and right).}
    \label{fig:space_dependency}
\end{figure}

Let us further introduce a space-dependency into the weight function of the DO differential operator. 
Then the system reads as 
\begin{equation}
    \int_0^3 \phi(\alpha, \Vec{x}) \tfC{\alpha} u(t,\Vec{x}) \Diff \alpha = \epsilon \Laplace u(t,\Vec{x}) + f(t,\Vec{x}), 
\end{equation}
subject to initial conditions 
\begin{alignat}{3}
u(0,\Vec{x}) 
= u_0(\Vec{x}), \quad 
u^{(1)}(0,\Vec{x})
\equiv 0, \quad 
u^{(2)}(0,\Vec{x})
\equiv 0, \quad 
\Vec{x} \in \Omega, 
\end{alignat}
and homogeneous Dirichlet boundary conditions. 
Note that we consider $\alpha \in [0,3]$, and thus, there is also an initial condition for the second time derivative of the solution. 
We employ the kernel compression method to point-wise approximate the space-dependent DO differential operator, whereby the weights and poles of the exponential sum approximations are space-dependent. 
For simplicity, we consider a fixed number of exponential terms $m=20$ for all integral kernels and for all $\Vec{x} \in \Omega$. 
The resulting system reads as  
\begin{alignat}{3}
    \sum_{i=1}^3 \sum_{j=1}^{m} v_{i,j}(t,\Vec{x}) & = \epsilon \Laplace v(t,\Vec{x}) + f(t,\Vec{x}), \\ 
    \frac{\mathrm{d}}{\mathrm{d}t} v_{i,j}(t,\Vec{x}) & = -\lambda_{i,j}(\Vec{x}) v_{i,j}(t,\Vec{x}) + w_{i,j}(\Vec{x}) v (t,\Vec{x}), \quad && j = 1, \ldots, m, \quad i = 1,2,3. 
\end{alignat}

We consider $\phi(\alpha, \Vec{x})$ to be a bump function with space-dependent center $\eta:\Omega\rightarrow[0,3]$ and compact support with respect to $\alpha$
\begin{equation}
    \phi(\alpha, \Vec{x}) = 
    \begin{cases}
        \tilde{c} \exp \left( \frac{1}{(\alpha-\eta(\Vec{x}))^2 - r^2}\right), & \alpha \in (\eta(\Vec{x})-r, \eta(\Vec{x})+r), \\ 
        0, & \text{otherwise,}
    \end{cases}
\end{equation}
whereby $\tilde{c}$ is a normalizing constant. 
We assume $\eta(\Vec{x})$ to take a discrete set of values and then compute the exponential sum approximation for all resulting integral kernels. 
The space-dependent weights and poles are taken as interpolating functions of these discrete values in the corresponding finite element space. 

We set $\Omega = (0,1)^2$, $T=3$, $\epsilon = 0.2$, $h =0.01$, $\Delta \Vec{x} = 1/512$, and consider multiple point sources of the form $u_0(\Vec{x}) = \sum_{l=1}^L \exp(-1000 \|\Vec{x} - \Vec{x_l}\|_2)$. 
Three different scenarios for the spatial dependency of $\eta(\Vec{x})$ are considered: a simple geometric setting, a random field, see \cite{duswald2024finite}, and a profile based on the Marmousi data set \cite{versteeg1994marmousi}. 
Note that, in the case of the Marmousi example, $\epsilon$ is not space-dependent, but a constant, highlighting the influence of the DO differential operator. 
For all three scenarios, the spatial dependency of $\eta(\Vec{x})$ as well as the location of the initial point sources are displayed in the left column of Figure \ref{fig:Marmousi}. 

The simple geometric setting is displayed in Figure \ref{fig:Geometric}. At the centers of the point sources, $\eta(\Vec{x})$ takes the values $0.5$, $1.17$, $1.83$, and $2.5$. 
Where $\eta(\Vec{x})<2$, the smaller $\eta(\Vec{x})$ the faster diffuse the point sources, while for $\eta(\Vec{x})=2.5$, the point expands and amplifies. 
In the case of $\eta(\Vec{x}) \in \{1.83, 2.5\}$, the point source develops into a wave profile. 
For both the random field (Figure \ref{fig:RandomField}) and the Marmousi (Figure \ref{fig:Marmousi}) scenarios, the solutions diffuse in areas where $\eta(\Vec{x})<2$ and amplify and propagate wave-like, in places where $\eta(\Vec{x})>2$. 

\subsection{Numerical example 4: Nonlinear elasticity with DO damping}

As a nonlinear DOFPDE, we consider the dynamical system describing a Neo-Hookean material with a DO damping term. 
The Neo-Hookean material law is given in terms of its strain energy \cite{bower2009applied}
\begin{equation}
    W(\tensor{J}) 
    = \frac{\mu}{2} (\det(\tensor{J})^{-2/3} \Trace(\tensor{J}\Transpose \tensor{J}) - 3) + \frac{\chi}{2}\left(det(\tensor{J})-1\right)^2, 
\end{equation}
whereby $\mu$ and $\chi$ denote the shear and bulk modulus, respectively. Accordingly, the 1st Piola--Kirchhoff stress tensor $\tensor{P}$ is given by 
\begin{equation}
    \tensor{P}(\tensor{J})
    = \mu \det(\tensor{J})^{-2/3} \tensor{J} 
    - \frac{\mu}{3} \det(\tensor{J})^{-2/3} \Trace(\tensor{J}\Transpose \tensor{J}) \tensor{J}^{-T} 
    + \chi \det(\tensor{J}) \tensor{J}^{-T} \left( \det(\tensor{J}) - 1\right). 
\end{equation}
Then, the dynamical system of the Neo-Hookean material with DO damping term reads as 
\begin{equation}
    \frac{\Diff^2}{\Diff t^2} \Vec{u}(t,\Vec{x}) 
    = \nabla \cdot \left(\tensor{P}\left( \nabla \Vec{u}(t,\Vec{x}) \right) 
    + \int_0^1 \phi(\alpha) \tfC{\alpha} \nabla \Vec{u}(t,\Vec{x}) \Diff \alpha\right) 
    + \Vec{f}(t,\Vec{x}), \quad \text{in } (0,T) \times \Omega, 
\end{equation}
where $\Vec{u}$ denotes the deformation field and $\Vec{f}$ an external force. The system is subject to initial conditions and mixed boundary conditions 
\begin{alignat}{3}
    \Vec{u}(0, \Vec{x}) 
    = \Vec{u_0}(\Vec{x}), \quad 
    \frac{\Diff}{\Diff t} \Vec{u}(0, \Vec{x}) & \equiv 0, \quad 
    && \Vec{x} \in \Omega, \\ 
    \tensor{P}\left(\nabla \Vec{u}(t, \Vec{x})\right) \Vec{n}(t, \Vec{x}) 
    & = \Vec{t_N}(t,\Vec{x}), \quad 
    && (t,\Vec{x}) \in (0,T) \times \Gamma_N, \\ 
    \Vec{u}(t, \Vec{x}) & \equiv \Vec{0}, \quad && (t,\Vec{x}) \in (0,T) \times \Gamma_D. 
\end{alignat}
After application of the kernel compression method, we obtain 
\begin{alignat}{3}
    \label{eq:NE_constraint}
    \frac{\Diff^2}{\Diff t^2} \Vec{v}(t,\Vec{x}) 
    & = \nabla \cdot \left(\tensor{P}\left( \nabla \Vec{v}(t,\Vec{x}) \right) 
    + \sum_{j=1}^m \nabla \Vec{v_j}(t,\Vec{x}) \Diff \alpha\right) 
    + \Vec{f}(t,\Vec{x}), \\ 
    \label{eq:NE_modes}
    \frac{\mathrm{d}}{\mathrm{d}t} \Vec{v_j} (t,\Vec{x}) 
    & = \lambda_{j} \Vec{v_j} (t,\Vec{x}) + w_j \Vec{v} (t,\Vec{x}), \quad 
    &&  j = 1, \ldots, m. 
\end{alignat}
The system is subject to initial conditions 
\begin{alignat}{3}
    \Vec{v}(0, \Vec{x}) 
    = \Vec{u_0}(\Vec{x}), \quad 
    \Vec{v}^{(1)}(0, \Vec{x}) 
    \equiv 0, \quad 
    \Vec{v_j}(0, \Vec{x}) 
    \equiv 0, \quad 
    \Vec{x} \in \Omega, \quad 
    j = 1, \ldots, m, 
\end{alignat}
and closed by the mixed boundary conditions 
\begin{alignat}{3}
    \tensor{P}\left(\nabla \Vec{v}(t, \Vec{x})\right) \Vec{n}(t, \Vec{x}) 
    & = \Vec{t_N}(t,\Vec{x}), \quad 
    && (t,\Vec{x}) \in (0,T) \times \Gamma_N, \\ 
    \Vec{v}(t, \Vec{x}) & \equiv \Vec{0}, \quad && (t,\Vec{x}) \in (0,T) \times \Gamma_D. 
\end{alignat}

\begin{figure}
    \centering
    \begin{subfigure}[b]{0.45\textwidth}
         \centering
         \includegraphics[width=\textwidth]{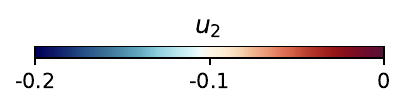}
    \end{subfigure}
    \hfill
    \begin{subfigure}[b]{0.45\textwidth}
         \centering
         \includegraphics[width=\textwidth]{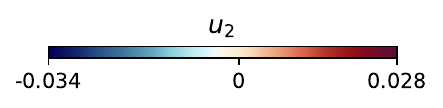}
    \end{subfigure}
    
    \begin{subfigure}[b]{0.49\textwidth}
         \centering
         \includegraphics[width=\textwidth, trim = {1cm 1cm 3cm 2cm}, clip]{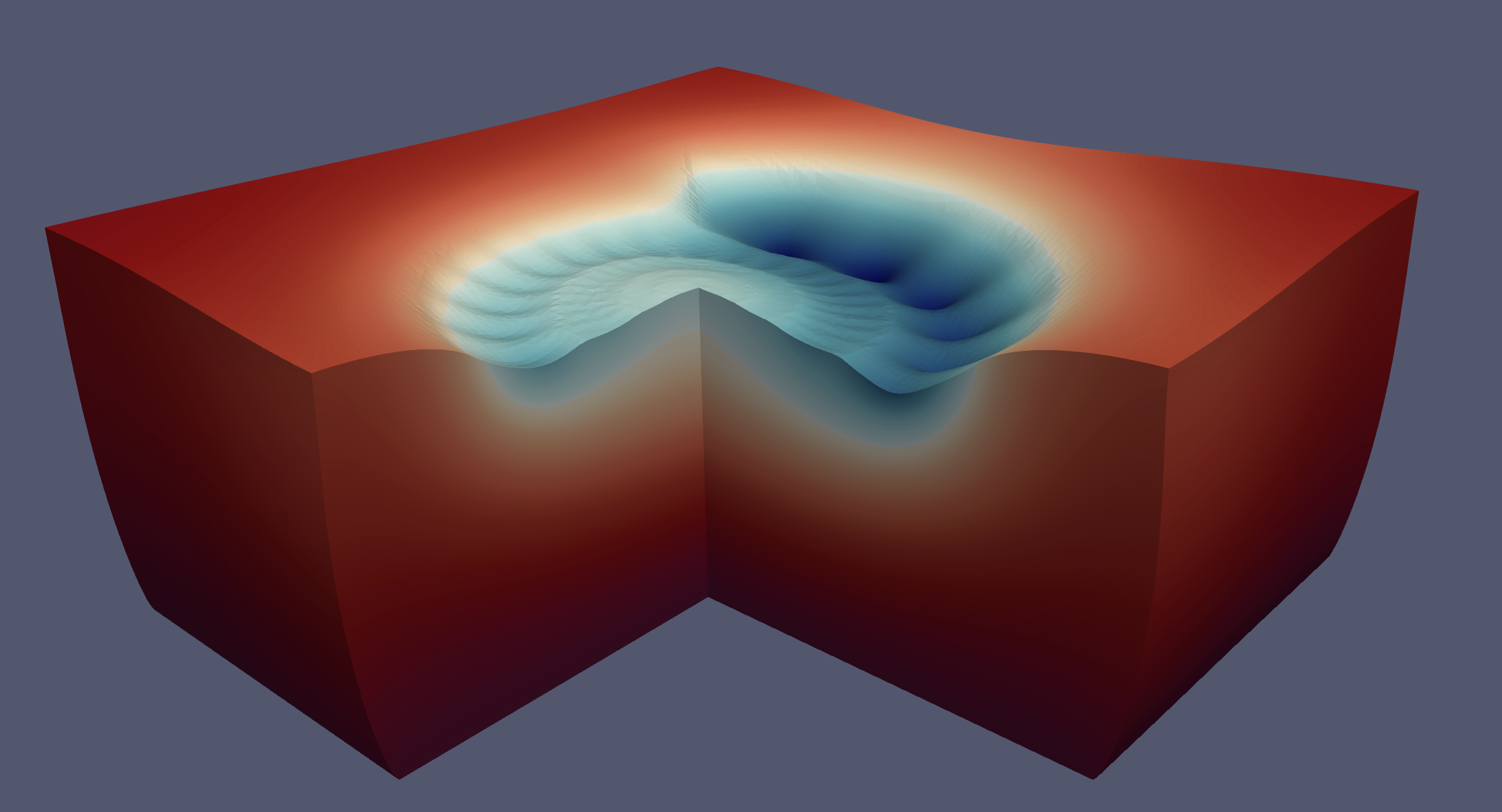}
         \caption{$\tau=0.1$, $t=2$}
    \end{subfigure}
    \hfill 
    \begin{subfigure}[b]{0.49\textwidth}
         \centering
         \includegraphics[width=\textwidth, trim = {1cm 1cm 3cm 2cm}, clip]{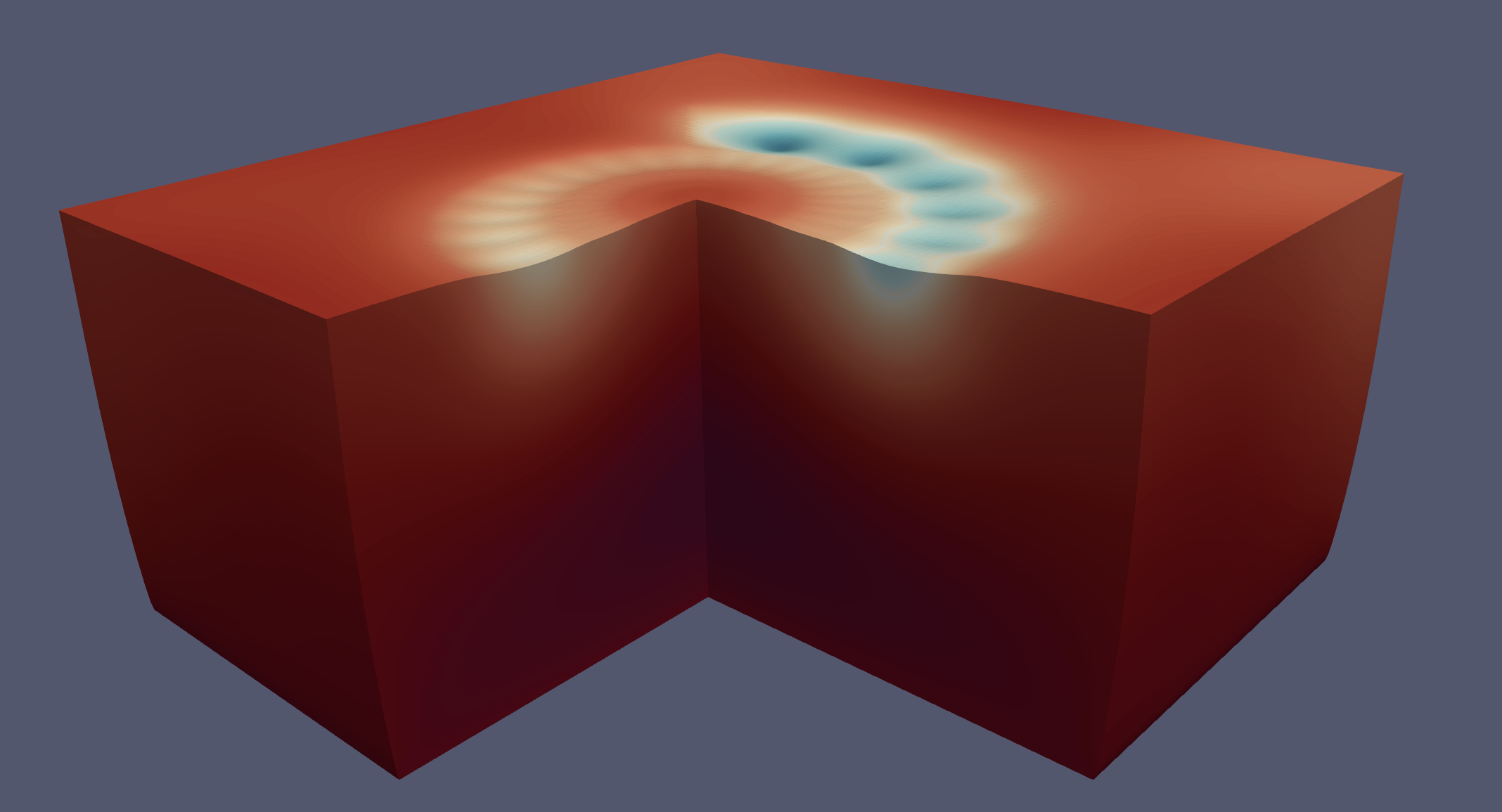}
         \caption{$\tau=0.1$, $t=2.5$}
    \end{subfigure}
    
    \begin{subfigure}[b]{0.49\textwidth}
         \centering
         \includegraphics[width=\textwidth, trim = {1cm 1cm 3cm 2cm}, clip]{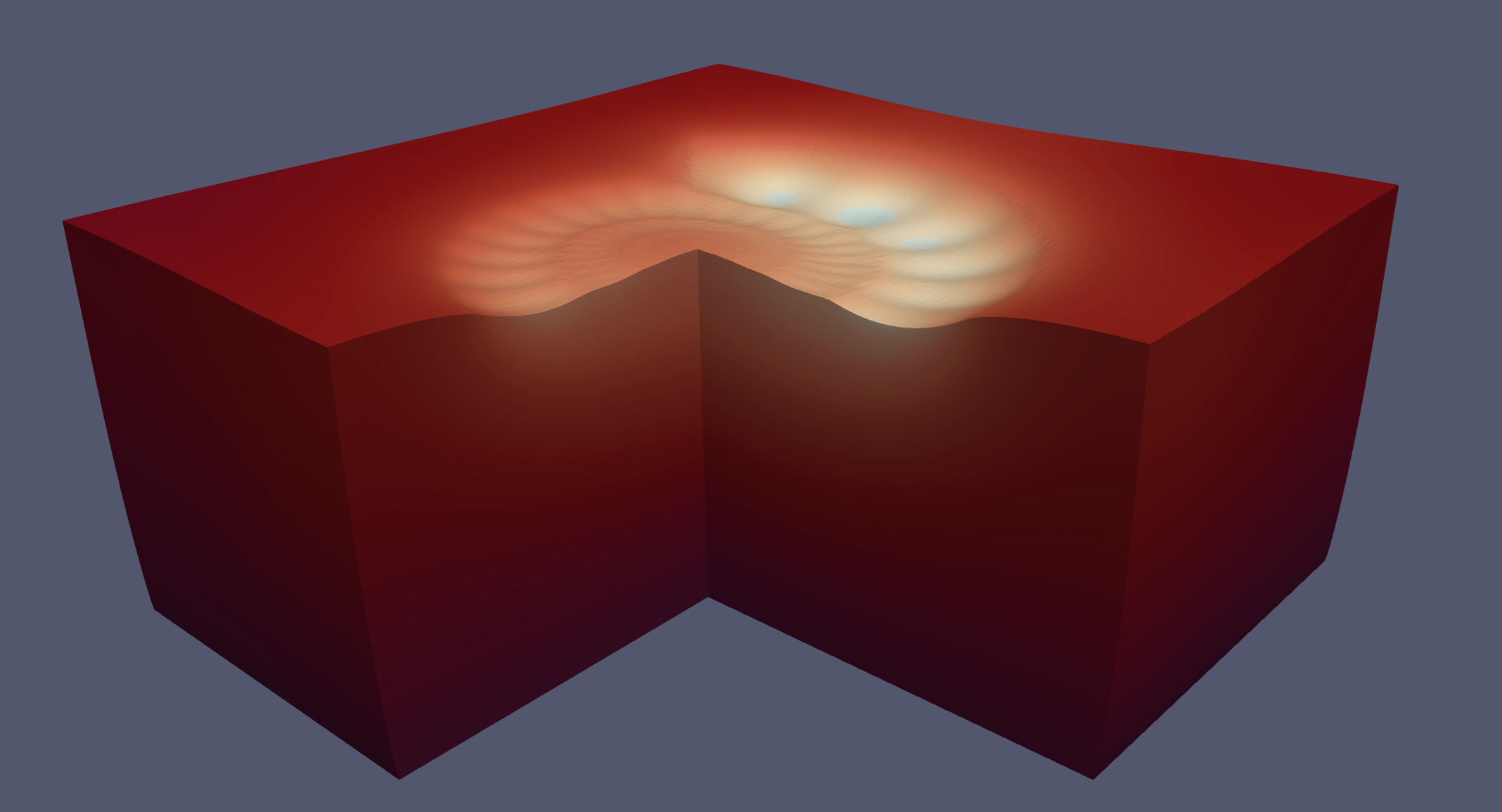}
         \caption{$\tau=2$, $t=2$}
    \end{subfigure}
    \hfill 
    \begin{subfigure}[b]{0.49\textwidth}
         \centering
         \includegraphics[width=\textwidth, trim = {1cm 1cm 3cm 2cm}, clip]{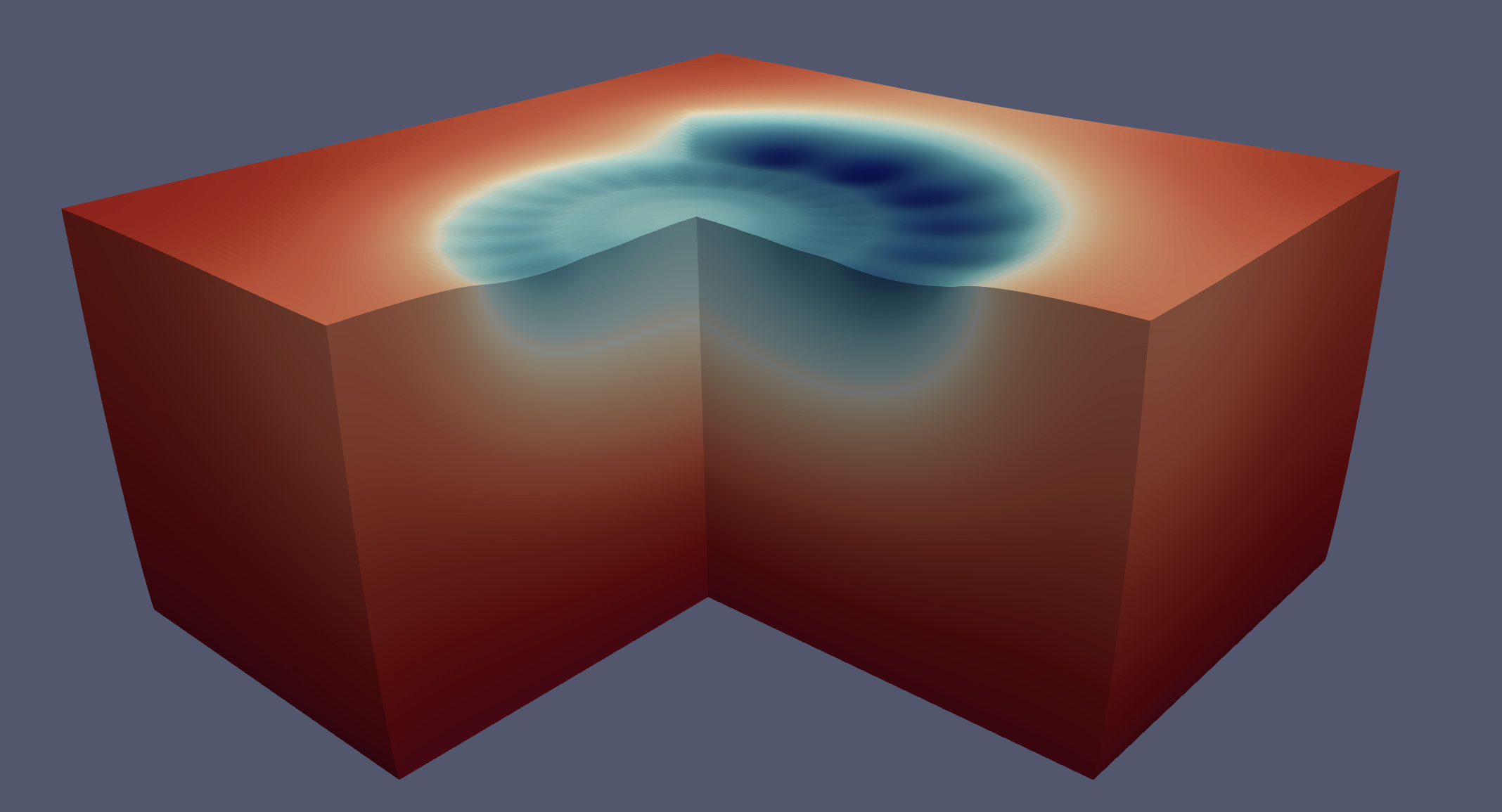}
         \caption{$\tau=2$, $t=2.5$}
    \end{subfigure}
    \caption{Deformation of a Neo-Hookean material with DO damping under load ($t=2$) and after release of the force ($t=2.5$) for two relaxation times $\tau = 0.1,2$.}
    \label{fig:compression}
\end{figure}

Let $\Omega = (0,1) \times (0,1) \times (0,0.5)$, $\mu=0.1$, $\chi=1$, $T=3$. We consider homogeneous Dirichlet boundary conditions and homogeneous Neumann boundary conditions at the bottom and the sides of the domain, respectively, and apply a constant external force in the negative $z$ direction in the shape of an ammonite. The force is multiplied with $t$, for $t \in [0,1]$, $1$ for $t \in [0,1]$, and the external force is set to $\Vec{0}$ otherwise. 
We consider DO differential operators with weight functions of the form $\tau^\alpha$, $\alpha \in [0,1]$, where $\tau > 0$ is the relaxation time of the damping \cite{atanackovic2005fractional}. 

We discretize \eqref{eq:NE_constraint}-\eqref{eq:NE_modes} using second-order finite elements and the implicit Euler method, and the resulting nonlinear system is linearized using Newton's method. 
The mesh consists of $500.000$ hexahedral elements, resulting in $12.241.503$ degrees of freedom, the timestep is set to $h=0.001$, and both integral kernels are approximated using $m=44$ exponential terms. 

The solution for $\tau\in\{0.1,2\}$ under loading ($t=2$) and after the release of the force ($t=2.5$) is displayed in Figure \ref{fig:compression}. 
In both cases, we observe substantial deformations in the $z$ direction under load and an elastic response after the release of the force. Compared to $\tau=0.1$, for $\tau=2$, the deformation under load is less pronounced and the relaxation is significantly slower. 

%% file: Files/Conclusion.tex
We presented a kernel compression method that offers a fast and memory-efficient framework enabling the numerical solution of high-dimensional, nonlinear DOFPDEs. 
By applying a quadrature formula, we transform the DO differential operators into a sum of fractional derivatives and replace the resulting linear combinations of the fractional integral kernels with exponential sums, approximating the nonlocal DO differential operator by a linear system of PDEs. 
We achieve exponential sum approximations of RL fractional integral kernels and the integral kernels arising from DO differential operators with double precision accuracy using less than $80$ and less than $92$ exponential terms, respectively. 
The moderate number of auxiliary variables facilitates the application of the kernel compression method to DOFPDEs without excessive storage requirements. 
After temporal discretization, due to the weakly coupled structure of the PDE system, static condensation reduces the PDE system to a single PDE that has to be solved at each time step, while the solution and the auxiliary variables are given by a linear update, not even requiring a matrix-vector multiplication. 
We applied our numerical scheme to solve the DO diffusion wave equation on a nontrivial three-dimensional domain with over $40$ million spatial degrees of freedom, approximate space-dependent DO differential operators, and to simulate a nonlinear, dynamical system describing the behavior of a Neo-Hookean material with DO damping.